\documentclass[12pt]{amsart}
\usepackage{amsmath,amssymb}
\usepackage{bbm}
\usepackage{graphicx,tikz,mathtools}
\usepackage[left=2.65cm,right=2.65cm,top=2.7cm,bottom=2.7cm]{geometry}
\usepackage[pdftex]{hyperref}
\usepackage{cite}
\usepackage{mathrsfs} 
\usepackage{caption}
\hypersetup{hidelinks}

\makeatletter
\@namedef{subjclassname@2020}{\textup{2020} Mathematics Subject Classification}
\makeatother

\newtheorem{theorem}{Theorem}[section]
\newtheorem{proposition}[theorem]{Proposition}

\newtheorem{lemma}[theorem]{Lemma}
\newtheorem*{conjecture*}{Conjecture}

\newcommand{\triple}[1]{{\left\vert\kern-0.25ex\left\vert\kern-0.25ex\left\vert #1
        \right\vert\kern-0.25ex\right\vert\kern-0.25ex\right\vert}}

\theoremstyle{definition}
\newtheorem{definition}[theorem]{Definition}

\newtheorem{remark}[theorem]{Remark}

\def\N{\mathbb{N}}
\def\R{\mathbb{R}}

\def\Z{\mathbb{Z}}

\newcommand\minus\backslash

\newcommand\lan\langle
\newcommand\ran\rangle

\renewcommand\leq\leqslant
\renewcommand\geq\geqslant

\numberwithin{equation}{section}

\allowdisplaybreaks

\begin{document}
 
\title[Schr\"odinger equation with Neumann conditions]{Nondegeneracy and uniqueness of ground states\\
of the Schr\"odinger equation\\ with Neumann conditions
\\ in a weakly nonlocal regime}

\author[S. Dipierro]{Serena Dipierro}
\address{ \vspace{-0.4cm}
\newline 
\textbf{{\small Serena Dipierro}}\nopagebreak[4]
\vspace{0.15cm}
\newline \indent Department of Mathematics and Statistics, University of Western Australia,
35 Stirling Highway, WA 6009 Crawley, Australia}
 \email{serena.dipierro@uwa.edu.au}
 
\author[M. Medina]{Maria Medina}
\address{ \vspace{-0.4cm}
\newline 
\textbf{{\small Maria Medina}}\nopagebreak[4]
\vspace{0.15cm}
\newline \indent  
Departamento de Matem\'aticas, Universidad Aut\'onoma de Madrid, 28049 Madrid, Spain}
 \email{maria.medina@uam.es}
 
\author[E. Valdinoci]{Enrico Valdinoci}
\address{ \vspace{-0.4cm}
\newline 
\textbf{{\small Enrico Valdinoci}}\nopagebreak[4]
\vspace{0.15cm}
\newline \indent Department of Mathematics and Statistics, University of Western Australia,
35 Stirling Highway, WA 6009 Crawley, Australia}
\email{enrico.valdinoci@uwa.edu.au}
 


\begin{abstract}  
We consider the nonlocal Schr\"odinger equation
driven by a fractional Laplacian of order~$2s$, with~$s$ close to~$1$.
Differently from the existing literature, which focused on either global problems in the whole of~$\R^n$ or Dirichlet conditions, we study here the case
of nonlocal Neumann conditions in the halfspace.

Moreover, we prove that ground states are nondegenerate, i.e. the only solutions of the linearized equation
are linear combinations of the tangential derivatives of the ground state.

Consequently, the ground state is unique, up to tangential translations and a change of sign.
\end{abstract}

\maketitle

\tableofcontents
 
\section{Introduction} \label{S.Introduction}

Nonlocal Schr\"odinger equations have a broad range of applications in quantum mechanics
and nonlinear analysis. Unlike the classical Schr\"odinger equation, where diffusion is governed by the Laplace operator, nonlocal models involve integro-differential operators, most notably the fractional Laplacian $(-\Delta)^s$, with $s\in(0,1)$, see e.g.~\cite{MR2932616}
and the appendix in~\cite{MR3393677} for detailed physical motivations
(see also~\cite{MR904142, MR2318846, MR3146377}
for similar equations arising in the description of boson stars).

The nonlocal nature of the fractional Laplacian introduces significant analytical challenges, since the
solutions exhibit a slow decay and the nonlocal character of the operator typically does not allow
the use of explicit solutions or {a} ``cut-and-paste'' method, see~\cite{MR1054115, MR3002595, MR3060858}.
Consequently, many techniques developed for local elliptic equations require substantial modifications.

One of the central topics in the study of nonlocal Schr\"odinger equations is the existence and qualitative behavior of localized solutions. In particular, concentration phenomena describe families of solutions that become increasingly localized around specific regions of the domain as a parameter tends to a limiting value. Such solutions are of considerable physical interest, as they represent semiclassical states where quantum particles concentrate and give rise to ``particles''.

So far, to the best of our knowledge, concentration phenomena for the nonlocal Schr\"odinger equations have been detected in the whole of~$\R^n$ (see~\cite{MR3059423,
MR3350616,
MR3494494,
MR4374612}),
and in bounded domains with Dirichlet external conditions (see~\cite{MR3393677}).

A fundamental ingredient in the analysis of concentration phenomena is the nondegeneracy of the ground states (i.e., roughly speaking, of solutions with minimal energy). Indeed, after rescaling around a concentration point, solutions of the singularly perturbed problem are expected to converge to a ground state, which thus {serves} as a convenient initial approximation for particle-like solutions. In this spirit, to construct genuine concentrating solutions, one typically regards the concentration problem as a small perturbation of the limiting one related to ground states and applies a bifurcation argument (such as a Lyapunov--Schmidt reduction). However, the success of this method, and in general of all methods based on {a} functional version of the ``implicit function theorem'',
heavily relies on the invertibility of the linearized operator around the ground state, modulo the natural kernel generated by the translational invariance of the equation. 

This property is usually called ``nondegeneracy''
and ensures also the local uniqueness and stability of the {ground state}, allowing the correction of approximate solutions into exact ones. Moreover, it yields precise asymptotic expansions of the reduced energy functional, from which the location and multiplicity of concentrating solutions can be often deduced according to the geometry of the potential.

Concentration problems in the whole of~$\R^n$ and in bounded domains with Dirichlet conditions (in which concentration can only take place at interior points) are thus intimately connected with the nondegeneracy of ground states in~$\R^n$. For the nonlocal Schr\"odinger equation,
this turned out to be a very complex topic of investigation, which has {now been resolved}, due to the contributions in~\cite{MR1111746, MR3070568, Fall-Valdinoci, MR3530361}.

The objective of this paper is to consider the nonlocal counterpart of this theory
for the halfspace in the presence of Neumann external data. For this, we will adopt a strategy introduced in~\cite{Fall-Valdinoci}, namely we will focus on the case in which the fractional parameter is close to~$1$, which enables the use of perturbative methods from the classical setting.

We stress that the case of Neumann external conditions is important both from a conceptual and an applied perspective. From the mathematical viewpoint, Neumann problems exhibit phenomena that are fundamentally different from those arising under Dirichlet boundary conditions. In particular, the absence of prescribed boundary values allows solutions to interact more strongly with the boundary, typically leading to new concentration patterns, also close to the boundary, and requiring different analytical techniques. Moreover, the associated variational framework possesses additional invariances and compactness issues that make the analysis substantially more delicate.

From the standpoint of applications, Neumann boundary conditions naturally model situations in which there is no flux across the boundary of the domain. Such conditions arise in quantum mechanics when the boundary acts as a perfectly reflecting barrier. Consequently, understanding concentration phenomena under Neumann boundary conditions is relevant for describing localized states in confined systems where particles or energy cannot escape.

The formulation of Neumann boundary conditions in the nonlocal setting is particularly significant,
since particles end up leaving a given domain in a finite time.
Consequently, a natural way to model confinement and ``no-flux conditions'' is to require that particles attempting to jump outside the domain are reflected back into it according to the same long-range jump distribution, see~\cite{MR3651008}.

Moreover, nontrivial solutions of the critical fractional Lane--Emden equation in the halfspace with these nonlocal Neumann conditions have recently been constructed in~\cite{Cinti-Talluri-Weth-2026,DelaTorre-Dipierro-Pistoia-Valdinoci-2026}.

{F}rom the technical standpoint, this kind of nonlocal condition is significantly harder than the classical one, since it involves an integral (rather than a differential) prescription,
placed in all the complement of a given domain (rather than along its boundary).
Also, while the classical Neumann condition can be easily ``localized'' (namely, if it is valid on~$\partial\Omega$,
then it is valid on~$B\cap\partial \Omega$, for every ball~$B$), the nonlocal Neumann condition intrinsically possesses a global character and localizations need to be carefully addressed.

Hence, in terms of applications, the study of nonlocal Schr\"odinger equations
in the halfspace subject to Neumann conditions is expected to be
substantially challenging, but also of notable importance, since it can pave the way
to the understanding of particle reflection phenomena and boundary concentration in the presence
of long-tail diffusion patterns.\medskip

The mathematical setting that we deal with in this paper goes as follows.
We consider the following nonlocal problem
\begin{equation}\label{problem}
\begin{cases}
(-\Delta)^su+u=|u|^{p-1}u\quad \mbox{ in }\R^n_+,\\
\mathcal{N}_su=0\quad \mbox{ in }\R^n_-,
\end{cases}
\end{equation}
where $s\in (0,1)$, $p\in(1,2^*_s-1)$, $2^*_s:=\frac{2n}{n-2s}$, $n\geq 3$, $\R^n_+:=\R^{n-1}\times(0,+\infty)$, $\R^n_-:=\R^{n-1}\times(-\infty,0)$ and 
$$\mathcal{N}_su(x):=c_{n,s}\int_{\R^n_+}\frac{u(x)-u(y)}{|x-y|^{n+2s}}\,dy$$
is the nonlocal Neumann condition given in~\cite{MR3651008}. {The coefficient~$c_{n,s}$ is the positive normalizing constant}\footnote{For the purposes of this paper, it suffices to know that~$c_{n,s}\sim1-s$
when~$s\sim1$.} of the fractional Laplacian, which, as customary, is defined as
$$(-\Delta)^s u(x):=c_{n,s}\lim_{\varepsilon\searrow0}
\int_{\R^n\setminus B_\varepsilon(x)}
\frac{u(x)-u(y)}{|x-y|^{n+2s}}\,dy,$$
see e.g.~\cite[formula~(3.1)]{MR2944369}.

In this paper, we assume that~$s$ is sufficiently close to~$1$ (depending on~$n$ and~$p$). This regime lies ``close to the classical case'' and will be referred to as ``weakly nonlocal''.

For the variational formulation, let us define
$$\|u\|_{s,+}:=\bigg(\frac{c_{n,s}}{2}\iint_{\R^{2n}\setminus(\R^n_-\times \R^n_-)}\frac{(u(x)-u(y))^2}{|x-y|^{n+2s}}\,dx\,dy+\int_{\R^n_+}u^2\,dx\bigg)^{\frac 12}$$
and
$$H^s_+:=\{u\mbox{ measurable such that }\|u\|_{s,+}<+\infty\}.$$
Given the variational nature of $\mathcal{N}_su$, solutions to \eqref{problem} correspond, up to multiplication by~$\nu_s^{\frac1{p-1}}$, to constrained critical points of 
the energy functional
$$\mathcal{J}_s(u,\nu):=\frac{1}{2}\|u\|_{s,+}^2-\frac{\nu}{p+1}\int_{\R^n_+}|u|^{p+1}\,dx,$$
and we define
\begin{equation}\label{MIN:GRO:ST}
\nu_s:=\inf_{{u\in H^s_+}\atop{\|u\|_{L^{p+1}(\R^n_+)}=1}}\|u\|_{s,+}^2.\end{equation}

\begin{definition}
We say that $u_s$ is a {\it ground state} solution of problem~\eqref{problem} 
if 
(up to multiplication by~$\nu_s^{\frac1{p-1}}$)
it {realizes} the previous infimum, i.e., if $\|u_s\|_{s,+}^2=\nu_s$ and $\|u_s\|_{L^{p+1}(\R^n_+)}=1$.
\end{definition}

The aim of this article is to prove nondegeneracy and uniqueness of ground state solutions\footnote{As customary in variational problems,
one can consider here solutions in the weak sense:
we say that $w_s$ is a weak solution of \eqref{problem_nu} if
$$(w_s,\varphi)_{s,+}=\nu_s \int_{\R^n_+}|w_s|^{p-1}w_s\varphi\,dx\qquad \mbox{ for every }\varphi \in H_+^s,$$
where\label{ihqkdncs:MAsxRT}
$$(u,v)_{s,+}:=\frac{c_{n,s}}{2}\iint_{\R^{2n}\setminus(\R^n_-\times \R^n_-)}\frac{(u(x)-u(y))(v(x)-v(y))}{|x-y|^{n+2s}}\,dx\,dy+\int_{\R^n_+}uv\,dx.$$}
of \eqref{problem}.

Consequently, $u_s$ solves 
\begin{equation}\label{problem_nu}
\begin{cases}
(-\Delta)^su_s+u_s=\nu_s|u_s|^{p-1}u_s\quad \mbox{ in }\R^n_+,\\
\mathcal{N}_su_s=0\quad \mbox{ in }\R^n_-,
\end{cases}
\end{equation}
and thus $v_s:=\nu_s^{\frac{1}{p-1}}u_s$ solves problem~\eqref{problem}. For convenience, we will work with the formulation given in~\eqref{problem_nu}. 

Moreover, for a positive ground state~$u_s$, the function $\partial_{x_i}u_s$, for any $i=1,\ldots, n-1$, is a solution of
\begin{equation}\label{lin_problem}
\begin{cases}
(-\Delta)^sw+w=\nu_spu_s^{p-1}w\quad \mbox{ in }\R^n_+,\\
\mathcal{N}_sw=0\quad \mbox{ in }\R^n_-,
\end{cases}
\end{equation}
that is, it is a solution to the linearized operator around $u_s$
{(see~\eqref{0qoejld23:TSBqP} in Proposition~\ref{decompw} for the proof).}
We will actually prove that, for $s$ close to 1, the kernel is spanned by~$\partial_{x_1}u_s,\dots,\partial_{x_{n-1}}u_s$.

\begin{theorem}[Nondegeneracy]\label{th:nondegeneracy}
Let $u_s$ be a {positive} normalized ground state solution of \eqref{problem_nu}. There exists $s_\star$ such that, for every $s\in (s_\star,1)$, if {$w\in H^s_+$ is a weak solution of~\eqref{lin_problem}} then
$$w\in \mbox{span}\{\partial_{x_1}u_s,\ldots, \partial_{x_{n-1}}u_s\}.$$
 \end{theorem}
 
 \begin{theorem}[Uniqueness]\label{th:uniqueness}
 There exists $s_\star$ such that, for every $s\in (s_\star,1)$, the ground state solution is unique, up to a translation in the first $n-1$ coordinates and multiplication by~$-1$.
 \end{theorem}
 
Alongside these main results, we provide a detailed description of the ground states in terms of existence, regularity and qualitative properties.
Specifically, we stress that the infimum in~{\eqref{MIN:GRO:ST}} is indeed attained:

\begin{theorem}[Existence]\label{EXIST:GRS}
There exists a ground state solution of problem~\eqref{problem}.
\end{theorem}
 
Additional information can be provided
on the ground states of problem~\eqref{problem}
(these properties also play a pivotal role in the proofs of the main
results on nondegeneracy and uniqueness).

\begin{theorem}[Fundamental properties of the ground states]\label{ADDITIO:PR}
Let~$u$ be a ground state solution of problem~\eqref{problem}.
Then, up to a sign change and a translation in~$\R^{n-1}$,
\begin{itemize}
\item For all~$x\in\R^n$,
\begin{equation}\label{REPOSIt} u(x)>0 .\end{equation}
\item There exists a function~$\overline u:[0,+\infty)\times\R\to\R$ such that,
for all~$(x',x_n)\in\R^{n-1}\times\R$ we have that
\begin{equation}\label{REARR:in} u(x',x_n)=\overline u(|x'|,x_n).\end{equation}
\item We have that 
\begin{equation}\label{OFCLA}
{\begin{gathered}
u\in L^\infty(\R^n)\cap C^\infty(\R^n_+),\\
u\in C^{\min\{2s,s+\frac12\}+\alpha}(\overline{\R^n_+}),
\end{gathered}}\end{equation}
for some~${\alpha\in(0,\frac12)}$.
\item For all~$\rho$, $x_n\in(0,+\infty)$,
\begin{equation}\label{eqRedde} \partial_\rho\overline u(\rho,x_n)<0.\end{equation}
\end{itemize}
\end{theorem}

We stress that Theorems~\ref{EXIST:GRS} and~\ref{ADDITIO:PR} hold true for all~$s\in(0,1)$.
The weakly nonlocal assumption that $s$ is close to~$1$ is used in this paper only
for Theorems~\ref{th:nondegeneracy} and~\ref{th:uniqueness}.

Let us now briefly discuss the methodological differences of this paper with respect to the existing literature. On the one hand, the overall perturbative strategy is shared with~\cite{Fall-Valdinoci}, since both this paper and that establish convergence as~$s\nearrow1$, pass to the classical linearized equation, and use classical nondegeneracy. 

On the other hand, the machinery needed for the exterior Neumann problem is significantly different. First of all, the variational framework differs since the energy functional considered here excludes interactions within~$\R^n_-\times\R^n_-$ and recovers the exterior values out of a Neumann-type (or Poisson-type) averaging formula. This is not merely a cosmetic difference, since the corresponding fractional problem here requires a bespoke boundary analysis and an additional mechanism is needed for preventing escape to infinity, relying on the comparison between halfspace and whole-space energies.

A second conceptual difference (also with respect to the classical case) deals with rearrangement methods, which must respect here the exterior Neumann condition. The strategy that we adopt consists in rearranging tangential slices and then defining an exterior rearranged function through the Neumann condition. This step, however, requires suitable integral estimates and a Riesz's rearrangement inequality applied to convolutions of kernels.

A third structural difference occurs in relation with compactness results, which, in this paper, need to address both a ``fixed'' and a ``translated'' halfspace configuration (the latter, possibly expanding to the whole of~$\R^n$). While the proofs of~\cite{Fall-Valdinoci} rely on global Fourier estimates and compactness of fully radial functions, these techniques are not immediately available in the halfspace (and tangential symmetry alone does not control transaltions in the normal direction).

A fourth important distinction is that~\cite{Fall-Valdinoci} employs the Caffarelli-Silvestre extension, while here we work directly with the integral energy, therefore the analytical mechanisms differ substantially. 

A fifth difference consists in how uniqueness is obtained. Specifically, while~\cite{Fall-Valdinoci} establishes uniform invertibility and a contraction argument, here we normalize the difference of two minimizers and reuse the limiting equation and weighted compactness to obtain a contradiction, thus providing a more direct, and somewhat elegant, uniqueness proof within the compactness framework already developed. 

Also, in terms of future avenues of investigation, the theorems
of this paper possess a strategic value, since they will play a pivotal role regarding the boundary concentration results for the nonlocal Schr\"odinger equations that will be presented in~\cite{FO:RTH}.

The rest of the paper is organized as follows.
We first prove the existence of ground
states and establish their basic properties: this
content is presented in Sections~\ref{EXIST:GRS:A}
and~\ref{EXIST:GRS:B}.

Then, in Section~\ref{LOC:MAXST1}
we localize the maximum of the ground state
uniformly as~$s$ approaches~$1$,
which plays a crucial role in the limit procedures needed
to prove the main results of nondegeneracy and uniqueness.

{Section~\ref{asymptotics} establishes convergence to the classical minimizer as~$s\nearrow1$. Sections~\ref{nondegeneracy:section} and~\ref{uniqueness:section} prove Theorems~\ref{th:nondegeneracy} and~\ref{th:uniqueness}, respectively.}

\section{Existence of the ground state and proof of Theorem~\ref{EXIST:GRS}}\label{EXIST:GRS:A}

We provide here the proof of Theorem~\ref{EXIST:GRS}. In a nutshell, the proof consists of the following steps:
\begin{itemize}
\item One picks a constrained minimizing sequence~$u_k$,
\item One proves that the sequence does not ``vanish'', namely it maintains a uniformly positive mass centered at some points~$y_k=(y_k',y_{k,n})\in\R^{n-1}\times(0,+\infty)$,
\item The coordinates~$y_{k,n}$ do not drift to infinity,
\item One can pass the minimizing sequence to the limit and obtain the desired ground state.
\end{itemize}
The method of proof is essentially by concentration-compactness, but we perform the calculations
in a self-contained way for the sake of completeness.

\begin{proof}[Proof of Theorem~\ref{EXIST:GRS}]
\noindent{\bf Step 1. Energy comparison and the minimizing sequence.}
We recall the fractional Sobolev inequality in~$H^s(\R^n)$, namely
\begin{equation}\label{SOBUTI} \|u\|_{L^{2^*_s}(\R^n)}\le C_n\left( c_{n,s} \iint_{\R^{2n}}\frac{(u(x)-u(y))^2}{|x-y|^{n+2s}}\,dx\,dy\right)^{\frac12},\end{equation}
where~$C_n\ge1$ is a constant depending only on~$n$ (we stress that~$C_n$ can be taken
independent of~$s$, thanks to the explicit value found in~\cite[Theorem~1.1]{MR2064421}).

We now obtain a similar estimate in~$\R^n_+$ by even reflection. Namely, we set
$$\widetilde u(x',x_n):=\begin{dcases}
u(x',x_n)&{\mbox{ if }}x_n\ge0,\\
u(x',-x_n)&{\mbox{ if }}x_n<0,
\end{dcases}$$and we observe that
\begin{align*}
&\|u\|_{L^{2^*_s}(\R^n_+)}
=\frac{ \|\widetilde u\|_{L^{2^*_s}(\R^n)}}{2^{\frac1{2^*_s}}}
\le \frac{C_n}{2^{\frac1{2^*_s}}}\left( c_{n,s} \iint_{\R^{2n}}\frac{(\widetilde u(x)-\widetilde u(y))^2}{|x-y|^{n+2s}}\,dx\,dy\right)^{\frac12}\\
&\quad=\frac{C_n\,\sqrt{c_{n,s}}}{2^{\frac1{2^*_s}-\frac12}}
\left(\begin{aligned}
&\iint_{\R^n_+\times\R^n_+}\frac{( u(x)- u(y))^2}{|x-y|^{n+2s}}\,dx\,dy\\
&+\iint_{\R^{n}_+\times\R^n_-}\frac{( u(x',x_n)- u(y',-y_n))^2}{|x-y|^{n+2s}}\,dx\,dy
\end{aligned}\right)^{\frac12}\\
&\quad=\frac{C_n\,\sqrt{c_{n,s}}}{2^{\frac1{2^*_s}-\frac12}}
\left(\begin{aligned}
&\iint_{\R^n_+\times\R^n_+}\frac{( u(x)- u(y))^2}{|x-y|^{n+2s}}\,dx\,dy\\
&+\iint_{\R^{n}_+\times\R^n_+}\frac{( u(x)- u(y))^2}{|(x'-y',x_n+y_n)|^{n+2s}}\,dx\,dy
\end{aligned}\right)^{\frac12}\\*
&\quad\le
2C_n\left( c_{n,s} \iint_{\R^n_+\times\R^n_+}\frac{(u(x)-u(y))^2}{|x-y|^{n+2s}}\,dx\,dy\right)^{\frac12}.
\end{align*}

Thus, setting~$\mu:=\frac{2^*_s-p-1}{2^*_s-2}$,
we observe that~$p+1=2\mu+2^*_s(1-\mu)$
and accordingly, using H\"older's inequality with exponents~$\frac1\mu$ and~$\frac1{1-\mu}$,
\begin{eqnarray*}&&
\|u\|_{L^{p+1}(\R^n_+)}\le
\|u\|_{L^{2}(\R^n_+)}^{\frac{2\mu}{p+1}}
\|u\|_{L^{2^*_s}(\R^n_+)}^{\frac{2^*_s(1-\mu)}{p+1}}\\&&\quad\le
(2C_n)^{\frac{2^*_s(1-\mu)}{p+1}} \|u\|_{L^{2}(\R^n_+)}^{\frac{2\mu}{p+1}}
\left( c_{n,s} \iint_{\R^n_+\times\R^n_+}\frac{(u(x)-u(y))^2}{|x-y|^{n+2s}}\,dx\,dy\right)^{\frac{2^*_s(1-\mu)}{2(p+1)}}.
\end{eqnarray*}

Next, let~$a:=\frac{p+1}{2\mu}$ and~$b:=\frac{p+1}{2^*_s(1-\mu)}$.
Since
$$\frac1a+\frac1b=\frac{2\mu}{p+1}+\frac{2^*_s(1-\mu)}{p+1}=1,
$$
we can use Young's inequality with exponents~$a$ and~$b$ to conclude that
\begin{eqnarray*}\frac{
\|u\|_{L^{p+1}(\R^n_+)}}{
(2C_n)^{\frac{2^*_s(1-\mu)}{p+1}}}&\le&\sqrt{ \|u\|_{L^{2}(\R^n_+)}^{\frac2{a}}
\left( c_{n,s} \iint_{\R^n_+\times\R^n_+}\frac{(u(x)-u(y))^2}{|x-y|^{n+2s}}\,dx\,dy\right)^{\frac1{b}}}
\\&\le&\sqrt{
\|u\|_{L^{2}(\R^n_+)}^2+c_{n,s} \iint_{\R^n_+\times\R^n_+}\frac{(u(x)-u(y))^2}{|x-y|^{n+2s}}\,dx\,dy}.
\end{eqnarray*}

Hence, since~$2^*_s(1-\mu)=\frac{n(p-1)}{2s}\le10n(p-1)$ as long as~$s\ge\frac1{20}$, we have that
$$ c^\star_{n,p}\|u\|_{L^{p+1}(\R^n_+)}^2\le
\|u\|_{L^{2}(\R^n_+)}^2+c_{n,s} \iint_{\R^n_+\times\R^n_+}\frac{(u(x)-u(y))^2}{|x-y|^{n+2s}}\,dx\,dy,$$
for some~$ c^\star_{n,p}>0$.

As a consequence, up to renaming constants,
\begin{equation}\label{NONZGSI}\nu_s\ge c^\star_{n,p}.
\end{equation}

In what follows, we will need to compare the infimum in~\eqref{MIN:GRO:ST}
with the corresponding quantity in the whole of~$\R^n$. More specifically, we define
\begin{equation}\begin{split}\label{def_nustar} \nu_s^\star&:=\inf_{{u\in H^s(\R^n)}\atop{\|u\|_{L^{p+1}(\R^n)}=1}}
\frac{c_{n,s}}{2}\iint_{\R^{2n}}\frac{(u(x)-u(y))^2}{|x-y|^{n+2s}}\,dx\,dy+\int_{\R^n}u^2\,dx
\\&=\inf_{u\in H^s(\R^n)\setminus\{0\}}\frac{\displaystyle
\frac{c_{n,s}}{2}\iint_{\R^{2n}}\frac{(u(x)-u(y))^2}{|x-y|^{n+2s}}\,dx\,dy+\int_{\R^n}u^2\,dx}{\displaystyle\|u\|_{L^{p+1}(\R^n)}^2}
.\end{split}\end{equation}
We know that the minimum in \eqref{def_nustar}
is attained (see~\cite{MR1054115, MR3002595, MR3060858}) by a radially symmetric function~$\omega_s$
with~$\|\omega_s\|_{L^{p+1}(\R^n)}^2=1$, such that
$$ \frac{\alpha_{n,s,p}}{1+|x|^{n+2s}}\le\omega_s(x)\le\frac{\beta_{n,s,p}}{1+|x|^{n+2s}},$$
for {suitable constants}~$\beta_{n,s,p}>\alpha_{n,s,p}>0$.

Given~$d>0$, to be taken conveniently large in what follows, possibly in dependence of~$n$, $s$, and~$p$, we let
$$\omega_{s,d}(x',x_n):=\omega_s(x',x_n-d).$$
We observe that if~$x_n<-d$ then
$$|x|^2=|x'|^2+x_n^2\ge |x'|^2+\frac{x_n^2}2+\frac{d^2}2\ge
\frac{d^2+|x|^2}2$$
and therefore,
using the change of variable~$X:=\frac{x}{d}$,
\begin{equation}\label{g.SK02w}\begin{split}
\|\omega_{s,d}\|_{L^{p+1}(\R^n_-)}^{p+1}&=
\int_{\R^{n-1}\times(-\infty,-d)} \omega_s^{p+1}(x)\,dx\\&
\le\beta_{n,s,p}^{p+1}\int_{\R^{n-1}\times(-\infty,-d)} \frac{dx}{|x|^{(n+2s)(p+1)}}\\&
\le C_{n,s,p}\int_{\R^{n-1}\times(-\infty,-d)} \frac{dx}{(d^2+|x|^2)^{\frac{(n+2s)(p+1)}2}}\\&=
\frac{ C_{n,s,p}}{d^{(n+2s)(p+1)-n}}\int_{\R^{n-1}\times(-\infty,-1)} \frac{dX}{(1+|X|^2)^{\frac{(n+2s)(p+1)}2}}\\&\le \frac{C_{n,s,p}}{d^{(n+2s)(p+1)-n}},
\end{split}
\end{equation}for some~$C_{n,s,p}\ge 1$, possibly varying from line to line.

On a similar note,
\begin{equation}\label{g.SK02w-dalb}\begin{split}&
\|\omega_{s,d}\|_{L^{2}(\R^n_-)}^{2}=
\int_{\R^{n-1}\times(-\infty,-d)} \omega_s^{2}(x)\,dx
\ge\int_{
B_d(0,\dots,0,-3d)
} \omega_s^2(x)\,dx\ge\frac{1}{C_{n,s,p}\,d^{n+4s}}.
\end{split}
\end{equation}

%%%%%Furthermore, if~$L_{n,s,p}:=5+3\left(\frac{4\beta_{n,s,p}}{\alpha_{n,s,p}}\right)^{\frac1{n+2s}}$, $
%%%%%e_n:=(0,\dots,0,1)$, $X\in B_1(-2e_n)$ and~$Y\in B_1(-L_{n,s,p} e_n)$, then~$|X|\le3$ and~$|Y|\ge4+3\left(\frac{4\beta_{n,s,p}}{\alpha_{n,s,p}}\right)^{\frac1{n+2s}}$, therefore
%%%%%\begin{eqnarray*}&& \omega_s(dX)-\omega_s(dY)\ge\frac{\alpha_{n,s,p}}{1+d^{n+2s}|X|^{n+2s}}-\frac{\beta_{n,s,p}}{1+d^{n+2s}|Y|^{n+2s}}\\&&\quad\ge
%%%%%\frac{\alpha_{n,s,p}}{3^{n+2s}d^{n+2s}}-\frac{2\beta_{n,s,p}}{d^{n+2s}|Y|^{n+2s}}
%%%%%\ge\frac{\alpha_{n,s,p}}{2\cdot 3^{n+2s}d^{n+2s}},\end{eqnarray*}
%%%%%as long as~$d$ is large enough.
%%%%%
%%%%%As a result,
%%%%%\begin{equation}\label{g.SK02w2}\begin{split}&
%%%%%\iint_{(\R^{n-1}\times(-\infty,-d))^2}\frac{(\omega_{s}(x)-\omega_{s}(y))^2}{|x-y|^{n+2s}}\,dx\,dy\\&\quad=
%%%%%d^{n-2s}\iint_{(\R^{n-1}\times(-\infty,-1))^2}\frac{(\omega_{s}(dX)-\omega_{s}(dY))^2}{|X-Y|^{n+2s}}\,dX\,dY
%%%%%\\&\quad\ge
%%%%%\frac{1}{C_{n,s,p}\,d^{n+6s}}\iint_{B_1(-2e_n)\times
%%%%%B_1(-L_{n,s,p} e_n)}\frac{dX\,dY}{|X-Y|^{n+2s}}\\&\quad
%%%%%\ge
%%%%%\frac{1}{C_{n,s,p}\,d^{n+6s}}.\end{split}
%%%%%\end{equation}

In addition,
setting~$t_d:=\int_{{\R^n_-}}\omega_{s,d}^{{p+1}}$, we see that
\begin{equation*}
 \|\omega_{s,d}\|_{L^{{p+1}}({\R^n_+})}^2=(1-t_d)^{2/{(p+1)}}\ge 1-t_d
\end{equation*}and therefore
\begin{equation}\label{EWNS:02oe-034.P}\begin{split}
\nu_s^\star\ge{\nu_s}(1-t_d)^{2/{(p+1)}}
             +\|\omega_{s,d}\|_{L^2({\R^n_-})}^2
 \ge{\nu_s}(1-t_d)+\|\omega_{s,d}\|_{L^2({\R^n_-})}^2>{\nu_s},
\end{split}\end{equation}
provided that~$d$ is chosen sufficiently large.

Now, to prove the existence of the desired ground state, up to scaling, it suffices to check that the infimum in~\eqref{MIN:GRO:ST} is attained. For this, we take a minimizing sequence~$u_k$, which is,
without loss of generality, nonnegative (up to replacing each function with its absolute value) and normalized such that~$\|u_k\|_{L^{p+1}(\R^n_+)}=1$.

Since the nonlocal Neumann condition decreases the interaction energy (see~\cite[Theorem~2.1]{MR4651677}),
without loss of generality we can suppose that, for all~$y\in\R^{n}_-$,
\begin{equation}\label{NSFQIFCELMI}
u_k(y)=\frac{\displaystyle\int_{\R^n_+}\frac{u_k(z)}{|y-z|^{n+2s}}\,dz}{\displaystyle\int_{\R^n_+}\frac{dz}{|y-z|^{n+2s}}}.
\end{equation}

\noindent{\bf Step 2. Nonvanishing.}
We claim that this sequence does not ``vanish'' as~$k\to+\infty$, namely that there exists~$\delta > 0$ such that
\begin{equation}\label{NOVAN}
\liminf_{k \to+\infty} \sup_{y \in{\R^n_+}} \int_{B_1(y) \cap \R^n_+} |u_k(x)|^{p+1} \, dx \geq \delta.\end{equation}
Assume for contradiction that this is not the case. Then, up to a subsequence,
\begin{equation}\label{NOVAN2}\lim_{k\to+\infty}\sup_{y\in{\R^n_+}} \int_{B_1(y) \cap \R^n_+} |u_k(x)|^{p+1} \, dx = 0.\end{equation}

We also set~$\theta:=\frac2{p+1}$ and we point out that~$\theta\in(0,1)$. In addition, we use the short notation
$$\sigma_k:=\sup_{y \in{\R^n_+}} \int_{B_1(y) \cap \R^n_+} |u_k(x)|^{p+1} \, dx.$$
Thus, we write~$\R^n_+$ as the non-overlapping union of cubes~$\{Q_j\}_{j\in\N}$ of diameter~$1$. In this way, using the fractional Sobolev inequality on each~$Q_j$ (see e.g.~\cite[Theorem~6.7]{MR2944369}), we find a constant~$C>0$, depending only on~$n$, $s$ and~$p$, possibly varying from line to line, such that
\begin{eqnarray*}
1&=&\int_{\R^n_+} |u_k(x)|^{p+1} \, dx\\&=&\sum_{j=0}^{+\infty}
\int_{Q_j} |u_k(x)|^{p+1} \, dx\\&\le&\sigma_k^{1-\theta}
\sum_{j=0}^{+\infty}
\left(\int_{Q_j} |u_k(x)|^{p+1} \, dx\right)^\theta
\\&
\le&
C\sigma_k^{1-\theta}
\sum_{j=0}^{+\infty}
\left(\iint_{Q_j\times Q_j}\frac{(u_k(x)-u_k(y))^2}{|x-y|^{n+2s}}\,dx\,dy+
\int_{Q_j} |u_k(x)|^{2} \, dx\right)^{\frac{(p+1)\theta}2}\\&=&
C\sigma_k^{1-\theta}
\sum_{j=0}^{+\infty}\left(
\iint_{Q_j\times Q_j}\frac{(u_k(x)-u_k(y))^2}{|x-y|^{n+2s}}\,dx\,dy+
\int_{Q_j} |u_k(x)|^{2} \, dx\right)\\&\le&
C\sigma_k^{1-\theta}
\sum_{j=0}^{+\infty}\left(
\iint_{\R^n\times Q_j}\frac{(u_k(x)-u_k(y))^2}{|x-y|^{n+2s}}\,dx\,dy+
\int_{Q_j} |u_k(x)|^{2} \, dx\right)\\&\le &C\sigma_k^{1-\theta}\|u_k\|_{s,+}^2.
\end{eqnarray*}
Therefore, taking the limit as~$k\to+\infty$ and recalling~\eqref{NOVAN2}, we deduce that~$1\le 0$, which is a contradiction. This completes the proof of~\eqref{NOVAN}.

As a consequence of~\eqref{NOVAN}, we can pick a sequence of points~$y_k=(y_k',y_{k,n})\in\R^{n-1}\times(0,+\infty)$ such that
\begin{equation}\label{ldmwj0jr43mpHS69mn8n7}
\int_{B_1(y_k) \cap \R^n_+} |u_k(x)|^{p+1} \, dx \geq \frac\delta2.\end{equation}
Also, up to replacing~$u_k(x)$ with~$u_k(x'+y_k',x_n)$, we can suppose that~$y_k'=0$.

\noindent{\bf Step 3. Control of the mass of the minimizing sequence and limits of monotone functions.}
We now claim that the points~$y_k$ do not drift to infinity, namely
\begin{equation}\label{NODRI}
{\mbox{if $y_k=(0,\dots,0,y_{k,n})$ satisfies~\eqref{ldmwj0jr43mpHS69mn8n7}
for some~$\delta>0$, then }}\,
\sup_{k\in\N} y_{k,n}<+\infty.
\end{equation}

Indeed, if not, up to a subsequence, we can suppose that
\begin{equation}\label{PSPMSykS-2p34}
\lim_{k\to+\infty} y_{k,n}=+\infty.
\end{equation}
To recenter the mass, we define~$v_k(x',x_n):=u_k(x',x_n+y_{k,n})$ and~$H_k:=\R^{n-1}\times(-\infty,-y_{k,n})$, and we observe that
\begin{equation}\label{SLsuS-4}
\nu_s= \lim_{k\to+\infty}\frac{c_{n,s}}{2} \iint_{\R^{2n} \setminus (H_k \times H_k)} \frac{(v_k(x)-v_k(y))^2}{|x-y|^{n+2s}} \, dx \, dy + \int_{\R^{n-1}\times(-y_{k,n},+\infty)} v_k^2(x) \, dx.
\end{equation}
%%In particular, for all balls~$B\subset\R^n$ and~$k$ large enough, we have that~$v_k$ is bounded in~$H^s(B)$. On this account, up to a subsequence, we have that~$v_k$ converges to some~$v_\infty$ weakly in~$H^s_{\rm{loc}}(\R^n)$,
%%strongly in~$L^{p+1}_{\rm{loc}}(\R^n)$, and a.e. in~$\R^n$ (see e.g.~\cite[Theorem 7.1]{MR2944369} for the relevant compact embedding).
%%
%%These observations and Fatou's Lemma give that
%%\begin{equation}\label{pk0q-wiLAKSv250X6-0}
%%\nu_s\ge\frac{c_{n,s}}{2} \iint_{\R^{2n} } \frac{(v_\infty(x)-v_\infty(y))^2}{|x-y|^{n+2s}} \, dx \, dy + \int_{\R^n} v_\infty^2(x) \, dx.
%%\end{equation}

%%%We point out that, by~\eqref{ldmwj0jr43mpHS69mn8n7}, 
%%%$$\|v_\infty\|_{L^{p+1}(B_1)}=\lim_{k\to+\infty}
%%%\|v_k\|_{L^{p+1}(B_1)}\ge\left(\frac{\delta}2\right)^{\frac1{p+1}}.$$
%%%
%%%%Moreover,
%%%%\begin{equation}\label{KPSAS-32rKS}
%%%%\begin{split}&\|v_\infty\|_{L^{p+1}(\R^n)}=\lim_{\varrho\to+\infty}
%%%%\|v_\infty\|_{L^{p+1}(B_\varrho)}=\lim_{\varrho\to+\infty}\lim_{k\to+\infty}
%%%%\|v_k\|_{L^{p+1}(B_\varrho)}\\&\;\le\lim_{\varrho\to+\infty}\lim_{k\to+\infty}
%%%%\|v_k\|_{L^{p+1}(\R^{n-1}\times(-y_{k,n},+\infty))}=\lim_{\varrho\to+\infty}\lim_{k\to+\infty}1=1.\end{split}\end{equation}

Now, given~$R>0$, we let
\begin{eqnarray*}&&Q_k(R) := \frac{c_{n,s}}{2}\iint_{(B_R \cap H_k^c) \times \R^n} \frac{(v_k(x)-v_k(y))^2}{|x-y|^{n+2s}} \, dx \, dy + \int_{B_R \cap H_k^c} v_k^2(x) \, dx\\&&\qquad\qquad\qquad\qquad+\int_{B_R \cap H_k^c} v_k^{p+1}(x) \, dx.
\end{eqnarray*}
Also, $Q_k$ is a nondecreasing function of~$R$. Moreover, $Q_k(R)$ is nonnegative and bounded uniformly in both~$k$ and~$R$, because,
in view of~\eqref{SLsuS-4} and taking~$k$ large enough,
\begin{align*}
Q_k(R)&\le\frac{c_{n,s}}{2} \iint_{\R^{2n} \setminus (H_k \times H_k)} \frac{(v_k(x)-v_k(y))^2}{|x-y|^{n+2s}} \, dx \, dy\\
&\quad+\int_{\R^{n-1}\times(-y_{k,n},+\infty)} v_k^2(x) \, dx+\|v_k\|_{L^{p+1}(H_k^c)}^{p+1}\\
&\le\nu_s+1+1.
\end{align*}
On this account, it follows from Helly's selection theorem that~$Q_k$
possesses a (not relabeled) pointwise convergent subsequence.
Thus, there exists a nondecreasing function~$Q_\infty:(0,+\infty)\to[0,\nu_s+1]$ such that,
for almost all $R>0$,
$$\lim_{k\to+\infty} Q_k(R) =Q_\infty(R).$$

Therefore, by the monotonicity and boundedness of~$Q_\infty$, we know that there exists~$L_\infty\in[0,\nu_s+1]$ such that
$$\lim_{R\to+\infty}Q_\infty(R)=L_\infty,$$
hence, for all~$m\in\N$, we can find~$\varrho_m>1$ such that, for all~$r$, $R\ge\varrho_m$, we have that~$|Q_\infty(r)-Q_\infty(R)|\le\frac1m$
(and, without loss of generality, we can think that~$\varrho_{m+1}>4\varrho_m$).

Owing to~\eqref{PSPMSykS-2p34}, we find~$K_m\in\N$ such that, for all~$k\ge K_m$, we have that~$y_{k,n}>10\varrho_{m+1}$
(and, without loss of generality, we can think that~$K_{m+1}>4K_m$).
Also, we pick~$R^\star_m\in\left[\varrho_m,\,\frac1{10}+\varrho_m\right]$ belonging to the full-measure
convergence set of~$Q_k$. Therefore, for all~$m\in\N$,
\begin{eqnarray*}&&
\frac1m\ge|Q_\infty(R^\star_m)-Q_\infty(R^\star_{m+1})|=
\lim_{k\to+\infty}|Q_k(R^\star_m)-Q_k(R^\star_{m+1})|.
\end{eqnarray*}
In this way, we find~$K^\star_m\ge K_m$ for which
$$\frac2m\ge|Q_{K^\star_m}(R^\star_m)-Q_{K^\star_m}(R^\star_{m+1})|,$$
and we
notice that~$B_{9R_{m+1}^\star}\subset H_{K_m^\star}^{c}$
and~$R_{m+1}^\star\ge2R_m^\star$.

That is, up to (not relabeled) subsequences,
we have a diverging sequence~$R_k$ for which
\begin{equation}\label{asdAaA}
B_{R_k}\subseteq B_{R_{k+1}}\subseteq H_k^c\end{equation} and
\begin{equation}\label{OMSLSUfnsaFRucufuR}\begin{split}
0&=\lim_{k\to+\infty}|Q_{k}(R_k)-Q_{k}(R_{k+1})|\\&=\lim_{k\to+\infty}\Bigg|
\frac{c_{n,s}}{2}\iint_{B_{R_k} \times \R^n} \frac{(v_k(x)-v_k(y))^2}{|x-y|^{n+2s}} \, dx \, dy + \int_{B_{R_k}} v_k^2(x) \, dx\\&\qquad\qquad\quad+ \int_{B_{R_k}} v_k^{p+1}(x) \, dx-
\frac{c_{n,s}}{2}\iint_{B_{R_{k+1}} \times \R^n} \frac{(v_k(x)-v_k(y))^2}{|x-y|^{n+2s}} \, dx \, dy \\&\qquad\qquad\quad- \int_{B_{R_{k+1}}} v_k^2(x) \, dx-
\int_{B_{R_{k+1}}} v_k^{p+1}(x) \, dx\Bigg|\\&=\lim_{k\to+\infty}
\frac{c_{n,s}}{2}\iint_{(B_{R_{k+1}}\setminus B_{R_k}) \times \R^n} \frac{(v_k(x)-v_k(y))^2}{|x-y|^{n+2s}} \, dx \, dy + \int_{B_{R_{k+1}}\setminus B_{R_k}} v_k^2(x) \, dx\\&\qquad\qquad\quad+\int_{B_{R_{k+1}}\setminus B_{R_k}} v_k^{p+1}(x) \, dx.
\end{split}\end{equation}

Now we pick~$\tau\in C^\infty_0(B_2,\,[0,1])$ with~$\tau=1$ in~$B_1$, and let~$\tau_k(x):=\tau\left(\frac{x}{R_k}\right)$. Let also~$v_{k,1}:=\tau_k v_k$ and~$v_{k,2}:=(1-\tau_k)v_k$.

By construction, $v_{k,1}\equiv0$ in~$\R^n\setminus B_{2R_k}$
and~$v_{k,2}\equiv0$ in~$B_{R_k}$. Therefore, we have that
$(v_{k,1}(x)-v_{k,1}(y))(v_{k,2}(x)-v_{k,2}(y))=0$ if~$(x,y)\in B_{R_k}\times B_{R_k}$ and if~$(x,y)\in(\R^n\setminus B_{2{R_k}})\times(\R^n\setminus B_{2{R_k}})$.

Consequently,
\begin{equation}\label{OMSLSUfnsaFRucufuR2}
\begin{split}
&\frac12\iint_{\R^{2n}\setminus (H_k \times H_k) }\frac{|v_{k,1}(x)-v_{k,1}(y)|\,|v_{k,2}(x)-v_{k,2}(y)|}{|x-y|^{n+2s}}\,dx\,dy\\&\quad\le
\iint_{(B_{R_k}\times(\R^n\setminus B_{2{R_k}}))}\frac{|v_{k,1}(x)-v_{k,1}(y)|\,|v_{k,2}(x)-v_{k,2}(y)|}{|x-y|^{n+2s}}\,dx\,dy\\&\qquad\quad
+\iint_{(B_{2{R_k}}\setminus B_{R_k})\times\R^n}\frac{|v_{k,1}(x)-v_{k,1}(y)|\,|v_{k,2}(x)-v_{k,2}(y)|}{|x-y|^{n+2s}}\,dx\,dy
\\&\quad\le
\iint_{B_{R_k}\times(\R^n\setminus B_{2R_k})}\frac{|v_{k,1}(x)|\,| v_{k,2}(y)|}{|x-y|^{n+2s}}\,dx\,dy\\&\qquad\quad+\sum_{j=1}^2\iint_{(B_{2{R_k}}\setminus B_{R_k})\times\R^n}\frac{|v_{k,j}(x)-v_{k,j}(y)|^2}{|x-y|^{n+2s}}\,dx\,dy.
\end{split}\end{equation}

In addition,
\begin{eqnarray*} |v_{k,j}(x)-v_{k,j}(y)|&\le&
|v_{k}(x)-v_{k}(y)|+|\tau_k(x)-\tau_k(y)|\,|v_k(x)|\\&\le&
|v_{k}(x)-v_{k}(y)|+
C\min\left\{ 1,\frac{|x-y|}{R_k}\right\}
\,|v_k(x)|,
\end{eqnarray*}where~$C>0$ may vary from line to line, leading to
\begin{equation*}\begin{split}&
\iint_{(B_{2{R_k}}\setminus B_{R_k})\times\R^n}\frac{|v_{k,j}(x)-v_{k,j}(y)|^2}{|x-y|^{n+2s}}\,dx\,dy
\\&\quad\le 2
\iint_{(B_{{R_{k+1}}}\setminus B_{R_k})\times\R^n}\frac{|v_{k}(x)-v_{k}(y)|^2}{|x-y|^{n+2s}}\,dx\,dy
+\frac{C}{R_k^{2s}}\int_{B_{{2R_{k}}}\setminus B_{R_k}} |v_{k}(x)|^2\,dx,
\end{split}
\end{equation*}
which is infinitesimal as~$k\to+\infty$, thanks to~\eqref{OMSLSUfnsaFRucufuR}.

{F}rom this and~\eqref{OMSLSUfnsaFRucufuR2} we arrive at
\begin{equation}\label{OKSL:SIK892SX}
\begin{split}&
\lim_{k\to+\infty}\iint_{\R^{2n}\setminus (H_k \times H_k) }\frac{|v_{k,1}(x)-v_{k,1}(y)|\,|v_{k,2}(x)-v_{k,2}(y)|}{|x-y|^{n+2s}}\,dx\,dy\\&\qquad\le2\lim_{k\to+\infty}
\iint_{B_{R_k}\times(\R^n\setminus B_{2R_k})}\frac{|v_{k}(x)|\,| v_{k}(y)|}{|x-y|^{n+2s}}\,dx\,dy
\\&\qquad\le2\lim_{k\to+\infty}
\iint_{B_{R_k}\times(\R^n\setminus B_{2R_k})}\frac{|v_{k}(x)|\,| v_{k}(y)-v_k(x)|}{|x-y|^{n+2s}}\,dx\,dy\\&\qquad\qquad\qquad+\iint_{B_{R_k}\times(\R^n\setminus B_{2R_k})}\frac{|v_{k}(x)|^2}{|x-y|^{n+2s}}\,dx\,dy.
\end{split}\end{equation}

Thus, since
\begin{eqnarray*}
&&\iint_{B_{R_k}\times(\R^n\setminus B_{2R_k})}\frac{|v_{k}(x)|^2}{|x-y|^{n+2s}}\,dx\,dy\le
\iint_{B_{R_k}\times(\R^n\setminus B_{R_k})}\frac{|v_{k}(x)|^2}{|z|^{n+2s}}\,dx\,dz\\&&\qquad
\le\frac{C}{R_k^{2s}}\int_{B_{R_k}}|v_{k}(x)|^2\,dx\le\frac{C(\nu_s+1)}{R_k^{2s}},
\end{eqnarray*}we deduce from~\eqref{OKSL:SIK892SX} that
\begin{equation*}
\begin{split}&
\lim_{k\to+\infty}\iint_{\R^{2n}\setminus (H_k \times H_k) }\frac{|v_{k,1}(x)-v_{k,1}(y)|\,|v_{k,2}(x)-v_{k,2}(y)|}{|x-y|^{n+2s}}\,dx\,dy
\\&\quad\le2\lim_{k\to+\infty}
\iint_{B_{R_k}\times(\R^n\setminus B_{2R_k})}\frac{|v_{k}(x)|\,| v_{k}(y)-v_k(x)|}{|x-y|^{n+2s}}\,dx\,dy\\
&\quad\le2\lim_{k\to+\infty}
\sqrt{\iint_{B_{R_k}\times(\R^n\setminus B_{2R_k})}\frac{|v_{k}(x)|^2}{|x-y|^{n+2s}}\,dx\,dy\;
\iint_{B_{R_k}\times(\R^n\setminus B_{2R_k})}\frac{| v_{k}(y)-v_k(x)|^2}{|x-y|^{n+2s}}\,dx\,dy}\\&\quad\le\lim_{k\to+\infty}\frac{C}{{R_k^s}}
\sqrt{\int_{B_{R_k}}|v_{k}(x)|^2\,dx\;
\iint_{B_{R_k}\times(\R^n\setminus B_{2R_k})}\frac{| v_{k}(y)-v_k(x)|^2}{|x-y|^{n+2s}}\,dx\,dy}\\&\quad\le{\lim_{k\to+\infty}\frac{C(\nu_s+1)}{R^s_k}},
\end{split}\end{equation*}which tends to zero as~$k\to+\infty$.

For this reason,
\begin{eqnarray*}
&& \iint_{\R^{2n} \setminus (H_k \times H_k)} \frac{(v_k(x)-v_k(y))^2}{|x-y|^{n+2s}} \, dx \, dy\\
&&\quad=
2\iint_{\R^{2n} \setminus (H_k \times H_k)} \frac{(v_{k,1}(x)-v_{k,1}(y))(v_{k,2}(x)-v_{k,2}(y))}{|x-y|^{n+2s}} \, dx \, dy\\&&\qquad\qquad+\sum_{j=1}^2\iint_{\R^{2n} \setminus (H_k \times H_k)} \frac{(v_{k,j}(x)-v_{k,j}(y))^2}{|x-y|^{n+2s}} \, dx \, dy\\&&\quad=\sum_{j=1}^2\iint_{\R^{2n} \setminus (H_k \times H_k)} \frac{(v_{k,j}(x)-v_{k,j}(y))^2}{|x-y|^{n+2s}} \, dx \, dy+o(1),
\end{eqnarray*}
as~$k\to+\infty$.

On top of that,
\begin{eqnarray*}
&& \int_{\R^{n-1}\times(-y_{k,n},+\infty)} v_k^2(x) \, dx\\&&\quad=
2\int_{\R^{n-1}\times(-y_{k,n},+\infty)} v_{k,1}(x)\,v_{k,2}(x) \, dx
+\sum_{j=1}^2
\int_{\R^{n-1}\times(-y_{k,n},+\infty)} v_{k,j}^2(x) \, dx
\\&&\quad=\sum_{j=1}^2
\int_{\R^{n-1}\times(-y_{k,n},+\infty)} v_{k,j}^2(x) \, dx+o(1)
,\end{eqnarray*}
as~$k\to+\infty$, owing to~{\eqref{OMSLSUfnsaFRucufuR}}.

Therefore, these observations and~\eqref{SLsuS-4} give that
\begin{equation}\label{TBRCmas3rf}
\begin{split}
\nu_s&= \lim_{k\to+\infty}\sum_{j=1}^2\left(\frac{c_{n,s}}{2} 
\iint_{\R^{2n} \setminus (H_k \times H_k)} \frac{(v_{k,j}(x)-v_{k,j}(y))^2}{|x-y|^{n+2s}} \, dx\,dy\right.\\
&\qquad \qquad\qquad
\left.+ \int_{\R^{n-1}\times(-y_{k,n},+\infty)} v_{k,j}^2(x) \, dx\right)
\\&\ge  \lim_{k\to+\infty}\sum_{j=1}^2 \nu_s\|v_{k,j}\|_{L^{p+1}(H_k^c)}^2
\end{split}\end{equation}
and therefore, recalling~\eqref{NONZGSI} and supposing, up to subsequences, that the limits exist,
\begin{equation}\label{TBEs0-pDELr0e} 1\ge\lim_{k\to+\infty}\sum_{j=1}^2 \|v_{k,j}\|_{L^{p+1}(H_k^c)}^2=
\alpha_1^2+\alpha_2^2,\end{equation}
where
$$\alpha_j:=\lim_{k\to+\infty}\|v_{k,j}\|_{L^{p+1}(H_k^c)}.$$

Moreover, for all~$a$, $b\ge0$, since~$p+1\ge2$ we have that
\begin{equation}\label{WEHi5oFI}\begin{split}&
(a+b)^{\frac{p+1}2}\ge a^{\frac{p+1}2}+b^{\frac{p+1}2},\\&
{\mbox{ with equality holding if and only if either~$a=0$ or~$b=0$.}}
\end{split}
\end{equation}
We thus deduce from~\eqref{TBEs0-pDELr0e} that
\begin{equation}\label{FEMafWka5t}\begin{split} 1&\ge
\alpha_1^{p+1}+\alpha_2^{p+1}\\&
=\lim_{k\to+\infty}\int_{H_k^c} \tau_k^{p+1}(x)\, |v_k(x)|^{p+1}\,dx
+\int_{H_k^c} (1-\tau_k(x))^{p+1}\, |v_k(x)|^{p+1}\,dx
.\end{split}\end{equation}

We also observe that
\begin{eqnarray*}&&
\int_{H_k^c} \big|\tau_k^{p+1}(x)-\tau_k(x)\big|\, |v_k(x)|^{p+1}\,dx
+\int_{H_k^c}\big| (1-\tau_k(x))^{p+1}-
(1-\tau_k(x))\big|\, |v_k(x)|^{p+1}\,dx\\&&\qquad
\le
2\int_{B_{2R_k}\setminus B_{R_k}} |v_k(x)|^{p+1}\,dx,
\end{eqnarray*}
which is infinitesimal as~$k\to+\infty$, thanks to~\eqref{OMSLSUfnsaFRucufuR}.

{F}rom this observation and~\eqref{FEMafWka5t}, we arrive at
\begin{equation*}\begin{split} &1\ge
\lim_{k\to+\infty}\int_{H_k^c} \tau_k(x)\, |v_k(x)|^{p+1}\,dx
+\int_{H_k^c} (1-\tau_k(x))\, |v_k(x)|^{p+1}\,dx\\&\qquad\qquad
=\lim_{k\to+\infty}\int_{H_k^c} |v_k(x)|^{p+1}\,dx=1
.\end{split}\end{equation*}
Consequently, recalling~\eqref{WEHi5oFI},
either~$\alpha_1=0$ or~$\alpha_2=0$.
However, in light of~\eqref{ldmwj0jr43mpHS69mn8n7}, we know that~$\alpha_1^{p+1}\ge\frac\delta2>0$,
and therefore~$\alpha_2=0$.

Hence, $\alpha_1=1$. We can thus reconsider~\eqref{TBRCmas3rf} and write that
\begin{equation}\label{cde4rfaki5PMSDNPKD}
\nu_s= \lim_{k\to+\infty}
\frac{c_{n,s}}{2} 
\iint_{\R^{2n} \setminus (H_k \times H_k)} \frac{(v_{k,1}(x)-v_{k,1}(y))^2}{|x-y|^{n+2s}} \, dx\,dy+ \int_{\R^{n-1}\times(-y_{k,n},+\infty)} v_{k,1}^2(x) \, dx.
\end{equation}

Moreover, $v_{k,1}\equiv0$ in~$H_k$, therefore~\eqref{cde4rfaki5PMSDNPKD} leads to
\begin{eqnarray*}
\nu_s&=& \lim_{k\to+\infty}
\frac{c_{n,s}}{2} 
\iint_{\R^{2n}} \frac{(v_{k,1}(x)-v_{k,1}(y))^2}{|x-y|^{n+2s}} \, dx\,dy+ \int_{\R^n} v_{k,1}^2(x) \, dx
\\&\ge&
\lim_{k\to+\infty}\|v_{k,1}\|_{L^{p+1}(H_k^c)}^2\,\nu^\star_s
\\&=&\nu_s^\star.
\end{eqnarray*}
This is in contradiction with~\eqref{EWNS:02oe-034.P} and the proof of~\eqref{NODRI} is completed.

\noindent{\bf Step 4. The limit of the minimizing sequence.}
Being a normalized minimizing sequence, we know that~$u_k$ is bounded in~$H^s_+$. On this account, up to a subsequence, we have that~$u_k$ converges to some~$u_\infty$ weakly in~$H^s_+$,
strongly in~$L^{p+1}(\R^n_+\cap B_\varrho)$ for all~$\varrho>0$, and a.e. in~$\R^n_+$ (see e.g.~\cite[Theorem~7.1]{MR2944369} for the relevant compact embedding). These observations and~\eqref{NSFQIFCELMI} also give that~$u_k\to u_\infty$
in~$\R^n$ (up to defining~$u_\infty$ as in~\eqref{NSFQIFCELMI} in~$\R^n_-$).

Thus, by Fatou's lemma,
\begin{equation}\label{ILPA:1}
\nu_s\ge\frac{c_{n,s}}{2} \iint_{\R^{2n}\setminus(\R^n_-\times\R^n_-) } \frac{(u_\infty(x)-u_\infty(y))^2}{|x-y|^{n+2s}} \, dx \, dy + \int_{\R^n_+} u_\infty^2(x) \, dx.
\end{equation}

Now we claim that
\begin{equation}\label{ILPA:2}
\|u_\infty\|_{L^{p+1}(\R^n_+)}=1.
\end{equation}
Indeed, if$$ \sigma:=\|u_\infty\|_{L^{p+1}(\R^n_+)}\quad{\mbox{ and }}\quad\tau:=\lim_{k\to+\infty}\|u_\infty-u_k\|_{L^{p+1}(\R^n_+)},$$ looking at the energy of~$u_\infty$ and~$u_\infty-u_k$, using the weak convergence of~$u_k$, we have that, as~$k\to+\infty$,
\begin{equation*}
 \|u_k\|_{s,+}^2=\|u_\infty\|_{s,+}^2+
                 \|u_k-u_\infty\|_{s,+}^2+o(1)
 \ge\nu_s\big(\|u_\infty\|_{L^{p+1}(\R^n_+)}^2+
                   \|u_k-u_\infty\|_{L^{p+1}(\R^n_+)}^2\big)+o(1).
\end{equation*}
Therefore, passing to the limit,
\begin{equation}\label{KS:oqld3}\nu_s\ge\lim_{k\to+\infty}
 \|u_k\|_{s,+}^2
 \ge\nu_s\big(\sigma^2+\tau^2\big).\end{equation}

In addition, by the Brezis--Lieb lemma,
\begin{equation*}\begin{split}&\sigma^{p+1}=
\int _{\R^n_+}|u_\infty(x)|^{p+1}\,dx =\lim _{k\to+\infty } \int _{\R^n_+}|u_k(x)|^{p+1}\,dx -\int _{\R^n_+}|u_\infty(x)-u_k(x)|^{p+1}\,dx
\\&\qquad\qquad\qquad=1-\lim _{k\to+\infty }\int _{\R^n_+}|u_\infty(x)-u_k(x)|^{p+1}\,dx=1-\tau^{p+1}.
\end{split}\end{equation*}

Therefore, combining this observation with~\eqref{KS:oqld3}
and recalling~\eqref{WEHi5oFI}, we conclude that
$$ 1=\sigma^{p+1}+\tau^{p+1}\le(\sigma^2+\tau^2)^{\frac{p+1}2}\le1,$$
and therefore either~$\sigma=0$ or~$\tau=0$.

However, by~\eqref{ldmwj0jr43mpHS69mn8n7}, we know that~$\sigma>0$. As a result, we deduce that~$\tau=0$,
and accordingly~$\sigma=1$.
This completes the proof of~\eqref{ILPA:2}.

In view of~\eqref{ILPA:1} and~\eqref{ILPA:2}, we see that~$u_\infty$ is the desired ground state.\end{proof}

\section{Additional properties of the ground state and proof of Theorem~\ref{ADDITIO:PR}}\label{EXIST:GRS:B}

In this section we provide the proof of Theorem~\ref{ADDITIO:PR}.

\begin{proof}[Proof of Theorem~\ref{ADDITIO:PR}]
The proof is quite technical, but it can be split into independent arguments.
Since the statements in Theorem~\ref{ADDITIO:PR}
are invariant under multiplication by a positive constant,
we will use the normalization~\eqref{problem_nu}, with
\begin{equation*}
\|u\|_{L^{p+1}(\R^n_+)}=1\qquad{\mbox{and}}\qquad
\|u\|_{s,+}^2=\nu_s.
\end{equation*}

\noindent{\bf Step 1. Preliminary positivity statements.}
We let~$u_+:=\max\{u,0\}$ and~$u_-:=\min\{u,0\}$. We observe that~$u=u_+ + u_-$ and
$$ (u_+(x)-u_+(y))^2+(u_-(x)-u_-(y))^2\le(u(x)-u(y))^2,$$
with equality holding if and only if~$u(x)u(y)\ge0$.

As a consequence,
\begin{equation*}\begin{aligned}
\nu_s&\ge
\frac{c_{n,s}}{2}\iint_{\R^{2n}\setminus(\R^n_-\times \R^n_-)}\frac{(u_+(x)-u_+(y))^2}{|x-y|^{n+2s}}\,dx\,dy\\
&\quad+\frac{c_{n,s}}{2}\iint_{\R^{2n}\setminus(\R^n_-\times \R^n_-)}\frac{(u_-(x)-u_-(y))^2}{|x-y|^{n+2s}}\,dx\,dy\\
&\quad+\int_{\R^n_+}u^2_+\,dx
+\int_{\R^n_+}u^2_-\,dx\\
&\ge\nu_s\big(\|{u_+}\|_{L^{p+1}(\R^n_+)}^2+\|{u_-}\|_{L^{p+1}(\R^n_+)}^2\big).
\end{aligned}\end{equation*}

Hence, by~\eqref{NONZGSI},
$$1\ge\|{u_+}\|_{L^{p+1}(\R^n_+)}^2+\|{u_-}\|_{L^{p+1}(\R^n_+)}^2.$$
We now recall~\eqref{WEHi5oFI}, used with~$a:=\|{u_+}\|_{L^{p+1}(\R^n_+)}^2$
and~$b:=\|{u_-}\|_{L^{p+1}(\R^n_+)}^2$, and we find that
$$1\ge\|{u_+}\|_{L^{p+1}(\R^n_+)}^{p+1}+\|{u_-}\|_{L^{p+1}(\R^n_+)}^{p+1}=\|u\|_{L^{p+1}(\R^n_+)}^{p+1}=1,
$$
showing that actually equality holds true in~\eqref{WEHi5oFI}, 
which says that necessarily either
\begin{equation*}
\|{u_+}\|_{L^{p+1}(\R^n_+)}^2=0
\qquad\mbox{or}\qquad
\|{u_-}\|_{L^{p+1}(\R^n_+)}^2=0.
\end{equation*}
{Thus, either~$u_+$ or~$u_-$ vanishes identically in~$\R^n_+$. After changing sign if necessary, we can assume that~$u\ge0$ in~$\R^n_+$. The Neumann condition in~\eqref{problem_nu} gives~$u\ge0$ in~$\R^n_-$ as well.}

\noindent{\bf Step 2. Boundedness, regularity and proof of~\eqref{REPOSIt} and~\eqref{OFCLA}}
{We take~$r$ such that~$p+1<r\le2_s^*$. An interpolation
inequality and a weighted Young inequality, together with~\eqref{SOBUTI}, give that, for all~$ w\in H^s_+$,
\begin{equation}\label{ground:sobolev}
\|w\|_{L^r(\R^n_+)}^2
\le\|w\|_{L^2(\R^n_+)}^2+
\|w\|_{L^{2_s^*}(\R^n_+)}^2
\le \widetilde{C}\|w\|_{s,+}^2,
\end{equation}
where~$\widetilde{C}$ depends only on~$n$.

Also, for~$N>0$, $\beta\ge1$ and~$a,b\ge0$, the Cauchy--Schwarz inequality gives
\begin{equation}\label{ground:truncated-powers}
\begin{aligned}
&\left|a\min\{a,N\}^{\beta-1}
-b\min\{b,N\}^{\beta-1}\right|^2\\
&\quad\le\frac{\beta^2}{2\beta-1}(a-b)
\left(a\min\{a,N\}^{2\beta-2}
-b\min\{b,N\}^{2\beta-2}\right).
\end{aligned}
\end{equation}
Moreover, the maps~$a\mapsto a\min\{a,N\}^{\beta-1}$ and~$a\mapsto a\min\{a,N\}^{2\beta-2}$ are Lipschitz and vanish at zero, so their compositions with~$u$ belong to~$H^s_+$. 

Therefore, recalling~\eqref{ground:truncated-powers}, testing~\eqref{problem_nu} with~$u\min\{u,N\}^{2\beta-2}$ and using H\"older's inequality with exponents~$\frac{p+1}{p-1}$ and~$\frac{p+1}{2}$, we see that
\begin{equation*}
\begin{split}
\|u\min\{u,N\}^{\beta-1}\|_{s,+}^2
&\le\beta\nu_s\int_{\R^n_+}u^{p+1}\min\{u,N\}^{2\beta-2}\,dx\\
&\le\beta\nu_s
\|u\|_{L^{p+1}(\R^n_+)}^{p-1}
\|u\min\{u,N\}^{\beta-1}\|_{L^{p+1}(\R^n_+)}^2\\
&=\beta\nu_s
\|u\min\{u,N\}^{\beta-1}\|_{L^{p+1}(\R^n_+)}^2.
\end{split}
\end{equation*}
From this, \eqref{ground:sobolev} and the monotone convergence theorem, we have that
\begin{equation}\label{ground:moser-step}
\|u\|_{L^{\beta r}(\R^n_+)}
\le(A\beta)^{1/(2\beta)}
\|u\|_{L^{\beta(p+1)}(\R^n_+)},
\end{equation}
where
\begin{equation}\label{bcjkegity4gw78yt348}
A:=\max\{1,\widetilde{C}\nu_s\}.\end{equation}

Now we let~$\kappa:=r/(p+1)>1$. Applying~\eqref{ground:moser-step} with~$\beta=\kappa^j$, for~$j\in\N$, gives
\begin{equation}\label{ground:moser-bound}
\|u\|_{L^\infty(\R^n_+)}
\le\prod_{j=0}^{+\infty}(A\kappa^j)^{1/(2\kappa^j)}
=A^{\frac{\kappa}{2(\kappa-1)}}
\kappa^{\frac{\kappa}{2(\kappa-1)^2}}<+\infty.
\end{equation}
This and the Neumann condition in~\eqref{problem_nu} give that
\begin{equation}\label{ground:moser-bound-2}
u\in L^\infty(\R^n).\end{equation}
as mentioned in~\eqref{OFCLA}.

We stress that, by the interior fractional regularity theory,
equation~\eqref{problem_nu} holds pointwise in~$\R^n_+$. Thus, if~$u(x_0)=0$ at an interior point~$x_0$, then
\begin{equation*}0=(-\Delta)^su(x_0)
=-c_{n,s}\int_{\R^n}\frac{u(y)}{|x_0-y|^{n+2s}}\,dy<0,
\end{equation*}which is
a contradiction. Thus~$u>0$ in~$\R^n_+$, and the Neumann condition gives~$u>0$ in~$\R^n_-$. 

In particular,
\begin{equation}\label{67PKS0ej2rL0akqj:r3}
u(x',x_n)>0
\quad\text{for }(x',x_n)\in\R^{n-1}\times(\R\setminus\{0\}).
\end{equation}
Moreover, since~$u$ is positive in~$\R^n_+$, interior bootstrap regularity for~\eqref{problem_nu} gives~$u\in C^\infty(\R^n_+)$
(as mentioned in~\eqref{OFCLA}).}

{Also, by~\eqref{ground:moser-bound-2} and~\cite[Theorem~4.1]{MR4630055}, we know that~$u\in C^\alpha(\overline{\R^n_+})$ (possibly after reducing~$\alpha\in(0,1)$). 

We also remark that the function~$t\mapsto\nu_st^p-t$ is Lipschitz on~$[0,\|u\|_{L^\infty(\R^n)}]$, therefore~$\nu_su^p-u\in C^\alpha(\overline{\R^n_+})$.
As a result, applying~\cite[Theorem~1.1, version~2]{2025arXiv251013340D} (with~$\alpha\in(0,1)$ reduced again if necessary), we obtain~$u\in C^{\min\{2s,s+\frac12\}+\alpha}(\overline{\R^n_+})$,
which completes the proof of~\eqref{OFCLA}. }

We now prove~\eqref{REPOSIt} by arguing for the sake of contradiction.
We use~\eqref{OFCLA} to define, by continuity,
$$u(x',0):=\lim_{x_n\searrow0}u(x',x_n)$$
and, by~\eqref{67PKS0ej2rL0akqj:r3}, we have that~$u\ge0$ in~$\R^n$
and, if~\eqref{REPOSIt} is violated, then there exists~$p'\in\R^{n-1}$ such that~$u(p',0)=0$. Also, using again~\eqref{OFCLA}, we have that there exists~$\rho>0$ such that~$\nu_su^{p-1}<1$ in~$\R^n_+\cap B_\rho(p',0)$.

Hence, using the fractional Hopf lemma (see e.g.~\cite[Lemma~1.2(i)]{MR3533199},
used here with~$c(x):=\nu_su^{p-1}(x)-1$ and~$\Omega:=\R^n_+\cap B_\rho(p',0)$), {we obtain}
\begin{equation}\label{Hopf bound} \liminf_{\tau\searrow0}\frac{u(p',\tau)}{\tau^s}>0.\end{equation}
{Additionally, since~$u(p',0)=0$, we infer from~\eqref{OFCLA} that
\begin{equation*}0\le\frac{u(p',\tau)}{\tau^s}
\le C\tau^{\min\{2s+\alpha,s+\frac12+\alpha,1\}-s},
\end{equation*}which is infinitesimal as~$\tau\searrow0$,
and this contradicts~\eqref{Hopf bound}. The proof
of~\eqref{REPOSIt} is thereby complete.}

\noindent{\bf Step 3. Rearrangement inequalities and proof of~\eqref{REARR:in}.}
We now consider, given~$x_n\in(0,+\infty)$, the Schwarz symmetric decreasing rearrangement~$u_\star(\cdot,x_n)$ of~$u(\cdot,x_n)$. This defines~$u_\star$ in~$\R^n_+$.
Also, for all~$y\in\R^n_-$, we set
\begin{equation*}
u_\star(y)=\frac{\displaystyle\int_{\R^n_+}\frac{u_\star(z)}{|y-z|^{n+2s}}\,dz}{\displaystyle\int_{\R^n_+}\frac{dz}{|y-z|^{n+2s}}}.
\end{equation*}

It is also useful to define, for all~$y\in\R^n_-$,
\begin{equation}\label{defAr7436983480} 
A(y):=\int_{\R^n_+}\frac{dz}{|y-z|^{n+2s}}=
\int_{\R^n_+}\frac{dZ}{|(Z',y_n-Z_n)|^{n+2s}}=\frac{C_{n,s}}{|y_n|^{2s}},
\end{equation}
where$$C_{n,s}:=\int_{\R^n_+}\frac{d\zeta}{|(\zeta',1+\zeta_n)|^{n+2s}},$$
and then, for all~$z\in\R^n$,
\begin{equation}\label{Pfy4i3truSDFGe763453278}
P(y,z):=
\frac{|y_n|^{2s}}{C_{n,s}\,|y-z|^{n+2s}}
=\frac{1}{A(y)|y-z|^{n+2s}}
.\end{equation}
In this way, for all~$y\in\R^n_-$,
\begin{equation*}
u_\star(y)=\frac{\displaystyle\int_{\R^n_+}\frac{u_\star(z)}{|y-z|^{n+2s}}\,dz}{A(y)}
=\int_{\R^n_+} u_\star(z)\,P(y,z)\,dz,
\end{equation*}
with
\begin{equation}\label{rearrangement:probability}
\int_{\R^n_+} P(y,z)\,dz=1.
\end{equation}
Similarly, for all~$y\in\R^n_-$,
\begin{equation*}
u(y)=\int_{\R^n_+} u(z)\,P(y,z)\,dz.
\end{equation*}

{Next, let
\begin{equation*}
\widetilde C_{n,s}:=\int_{\R^{n-1}}
\frac{d\zeta'}{(|\zeta'|^2+1)^{(n+2s)/2}}.
\end{equation*}
Also, for~$x_n\ne y_n$, the change of variables~$\zeta'=(y'-x')/|x_n-y_n|$ gives
\begin{equation}\label{rearrangement:tangential-integral}
\int_{\R^{n-1}}
\frac{dy'}{(|x'-y'|^2+|x_n-y_n|^2)^{(n+2s)/2}}
=\frac{\widetilde C_{n,s}}{|x_n-y_n|^{1+2s}}.
\end{equation}

Moreover, tangential rearrangement preserves the~$L^2$ norm on almost every slice, namely, for almost every~$x_n>0$,
\begin{equation}\label{rearrangement:slice-L2}
\int_{\R^{n-1}}u_\star^2(x',x_n)\,dx'
=\int_{\R^{n-1}}u^2(x',x_n)\,dx'.
\end{equation}
Therefore,~$\|u_\star\|_{L^2(\R^n_+)}=\|u\|_{L^2(\R^n_+)}<+\infty$.

Next, pick~$0<\epsilon<R$ and take either~$U=u$ or~$U=u_\star$.
We have that, for every~$y\in\R^n_-$,
\begin{equation}\label{rearrangement:neumann-average}
\int_{\R^n_+}\frac{U(x)}{|x-y|^{n+2s}}\,dx=A(y)U(y).
\end{equation}
Therefore, \eqref{Pfy4i3truSDFGe763453278}, \eqref{rearrangement:probability} and~\eqref{rearrangement:neumann-average}
give, for every~$y\in\R^n_-$,
\begin{equation}\label{rearrangement:jensen}
\begin{split}
A(y)U^2(y)
&=\frac1{A(y)}
\left(\int_{\R^n_+}\frac{U(x)}{|x-y|^{n+2s}}\,dx\right)^2\\
&=A(y)\left(\int_{\R^n_+}U(x)P(y,x)\,dx\right)^2\\
&\le A(y)\int_{\R^n_+}U^2(x)P(y,x)\,dx\\
&=\int_{\R^n_+}\frac{U^2(x)}{|x-y|^{n+2s}}\,dx,
\end{split}
\end{equation}
where Jensen's inequality has also been used.

Moreover, using~\eqref{rearrangement:tangential-integral} yields
\begin{equation}\label{rearrangement:strip-bound}
\begin{split}
&\iint_{\R^n_+\times\{-R<y_n<-\epsilon\}}
\frac{U^2(x)}{|x-y|^{n+2s}}\,dx\,dy
=\widetilde C_{n,s}\int_{\R^n_+}U^2(x)
\left(\int_{-R}^{-\epsilon}\frac{dy_n}{(x_n-y_n)^{1+2s}}\right)\,dx\\
&\qquad\qquad\le\frac{\widetilde C_{n,s}}{2s\epsilon^{2s}}
\|U\|_{L^2(\R^n_+)}^2<+\infty.
\end{split}
\end{equation}
As a consequence of~\eqref{rearrangement:jensen} and~\eqref{rearrangement:strip-bound}, we obtain that
\begin{equation}\label{rearrangement:squares-finite}
\begin{split}
&\iint_{\R^n_+\times\{-R<y_n<-\epsilon\}}
\frac{U^2(y)}{|x-y|^{n+2s}}\,dx\,dy
=\iint_{\R^n_+\times\{-R<y_n<-\epsilon\}}
U^2(y)A(y)P(y,x) \,dx\,dy
\\&\qquad
=\int_{\{-R<y_n<-\epsilon\}}A(y)U^2(y)\,dy
\le\iint_{\R^n_+\times\{-R<y_n<-\epsilon\}}
\frac{U^2(x)}{|x-y|^{n+2s}}\,dx\,dy<+\infty.
\end{split}
\end{equation}

Furthermore, by~\eqref{REPOSIt} and the definition of~$u_\star$ we see that~$U\ge0$.
Therefore, Tonelli's theorem,~\eqref{rearrangement:neumann-average} and~\eqref{rearrangement:squares-finite} give
\begin{equation}\label{rearrangement:product-identity}
\begin{aligned}
0&\le\iint_{\R^n_+\times\{-R<y_n<-\epsilon\}}
\frac{U(x)U(y)}{|x-y|^{n+2s}}\,dx\,dy\\
&=\int_{\{-R<y_n<-\epsilon\}}U(y)
\left(\int_{\R^n_+}\frac{U(x)}{|x-y|^{n+2s}}\,dx\right)\,dy\\
&=\int_{\{-R<y_n<-\epsilon\}}A(y)U^2(y)\,dy<+\infty.
\end{aligned}
\end{equation}
In light of~\eqref{rearrangement:squares-finite} and~\eqref{rearrangement:product-identity} we can thereby write
\begin{equation}\label{oslmwq03.aA}
\begin{split}
&\iint_{\R^n_+\times\{-R<y_n<-\epsilon\}}
\frac{(U(x)-U(y))^2}{|x-y|^{n+2s}}\,dx\,dy\\
&\quad=\iint_{\R^n_+\times\{-R<y_n<-\epsilon\}}
\frac{U^2(x)+U^2(y)-2U(x)U(y)}{|x-y|^{n+2s}}\,dx\,dy\\
&\quad=\iint_{\R^n_+\times\{-R<y_n<-\epsilon\}}
\frac{U^2(x)}{|x-y|^{n+2s}}\,dx\,dy
-\int_{\{-R<y_n<-\epsilon\}}U^2(y)A(y)\,dy.
\end{split}
\end{equation}

Moreover,~\eqref{rearrangement:tangential-integral} and~\eqref{rearrangement:slice-L2} give
\begin{equation}\label{nshj3ortg:pa123LSP:01peokf}
\begin{split}
&\iint_{\R^n_+\times\{-R<y_n<-\epsilon\}}
\frac{u^2(x)-u_\star^2(x)}{|x-y|^{n+2s}}\,dx\,dy\\
&\quad=\widetilde C_{n,s}\int_0^{+\infty}\int_{-R}^{-\epsilon}
\frac{\displaystyle\int_{\R^{n-1}}
\big(u^2(x',x_n)-u_\star^2(x',x_n)\big)\,dx'}
{(x_n-y_n)^{1+2s}}\,dy_n\,dx_n=0.
\end{split}
\end{equation}

Next, to compare the tangential~$L^2$ norms of~$u$ and~$u_\star$ in~$\R^n_-$, let~$y_n<0$ and define
\begin{equation*}\begin{split}
f_{y_n,x_n}(\xi')&:=\frac1{(|\xi'|^2+(x_n-y_n)^2)^{(n+2s)/2}}\\
{\mbox{and }}\qquad
g_{y_n,z_n}(\xi')&:=\frac1{(|\xi'|^2+(z_n-y_n)^2)^{(n+2s)/2}}.
\end{split}
\end{equation*}
We observe that the convolution~$f_{y_n,x_n}*g_{y_n,z_n}$ is radial decreasing, by~\cite[page~171]{MR385456}. Also, \eqref{rearrangement:neumann-average} and Tonelli's theorem give, for either~$U=u$ or~$U=u_\star$,
\begin{equation}\label{AOKSM.Sw3r4}
\begin{split}
&\int_{\R^{n-1}}U^2(y',y_n)\,dy'\\
&=\frac{|y_n|^{4s}}{C_{n,s}^2}
\int_0^{+\infty}\int_0^{+\infty}
\iint_{\R^{n-1}\times\R^{n-1}}
U(x',x_n)U(z',z_n)\\
&\qquad\cdot(f_{y_n,x_n}*g_{y_n,z_n})(x'-z')
\,dx'\,dz'\,dx_n\,dz_n.
\end{split}
\end{equation}Thus, applying
Riesz's rearrangement inequality
in~\cite[Theorem~3.7]{MR1817225} to the tangential integrals in~\eqref{AOKSM.Sw3r4}, we infer that
\begin{equation}\label{PS:023o0e2kr3pf:9o2rj}
\int_{\R^{n-1}}u^2(y',y_n)\,dy'
\le\int_{\R^{n-1}}u_\star^2(y',y_n) \,dy'.
\end{equation}

Since~$A(y)=C_{n,s}|y_n|^{-2s}$ is independent of~$y'$, equations~\eqref{oslmwq03.aA},~\eqref{nshj3ortg:pa123LSP:01peokf} and~\eqref{PS:023o0e2kr3pf:9o2rj} imply that
\begin{equation*}
\iint_{\R^n_+\times\{-R<y_n<-\epsilon\}}
\frac{(u(x)-u(y))^2}{|x-y|^{n+2s}}\,dx\,dy
\ge
\iint_{\R^n_+\times\{-R<y_n<-\epsilon\}}
\frac{(u_\star(x)-u_\star(y))^2}{|x-y|^{n+2s}}\,dx\,dy.
\end{equation*}
Consequently, letting~$R\to+\infty$ and~$\epsilon\searrow0$, the monotone convergence theorem yields
\begin{equation}\label{nshj3ortg:pa123LSP:01peokf.l}
\iint_{\R^n_+\times\R^n_-}
\frac{(u(x)-u(y))^2}{|x-y|^{n+2s}}\,dx\,dy
\ge
\iint_{\R^n_+\times\R^n_-}
\frac{(u_\star(x)-u_\star(y))^2}{|x-y|^{n+2s}}\,dx\,dy.
\end{equation}}

Also, for almost every pair~$x_n,y_n>0$ with~$x_n\ne y_n$, applying~\eqref{rearrangement:tangential-integral} and~\eqref{rearrangement:slice-L2} gives
\begin{equation}\label{rearrangement:regional-squares}
\iint_{\R^{n-1}\times\R^{n-1}}
\frac{u^2(x',x_n)-u_\star^2(x',x_n)}
{(|x'-y'|^2+|x_n-y_n|^2)^{(n+2s)/2}}\,dx'\,dy'=0.
\end{equation}

As a result, using~\eqref{rearrangement:regional-squares} also with~$x$ and~$y$ interchanged, together with the expansion
%by virtue of the fractional {P}\'{o}lya-{S}zeg\"{o} Inequality
%(see e.g.~\cite[Corollary~2.3]{MR1002633}
%or~\cite[Theorem A.1]{MR2469027}),
\begin{eqnarray*}&&
\iint_{\R^{n-1}\times\R^{n-1}}\frac{(U(x',x_n)-U(y',y_n))^2}{|x-y|^{n+2s}}\,dx'\,dy'\\&&\qquad=
\iint_{\R^{n-1}\times\R^{n-1}}\frac{U^2(x',x_n)+U^2(y',y_n)-2U(x',x_n)\,U(y',y_n)}{|x-y|^{n+2s}}\,dx'\,dy',
\end{eqnarray*}
we conclude that
\begin{equation}\label{rearrangement:regional}\begin{split}&
\iint_{\R^{n-1}\times\R^{n-1}}\frac{(u(x',x_n)-u(y',y_n))^2}{|x-y|^{n+2s}}\,dx'\,dy'\\
&\qquad-\iint_{\R^{n-1}\times\R^{n-1}}\frac{(u_\star(x',x_n)-u_\star(y',y_n))^2}{|x-y|^{n+2s}}\,dx'\,dy'\\&\qquad=
2\iint_{\R^{n-1}\times\R^{n-1}}\frac{u_\star(x',x_n)\,u_\star(y',y_n)-u(x',x_n)\,u(y',y_n)}{|x-y|^{n+2s}}\,dx'\,dy'\ge0,
\end{split}\end{equation}where Riesz's rearrangement inequality has been used again in the last step.

{Therefore, integrating~\eqref{rearrangement:regional} in~$x_n$, $y_n>0$ and using~\eqref{nshj3ortg:pa123LSP:01peokf.l}, we conclude that~$\|u_\star\|_{s,+}\le\|u\|_{s,+}$. In fact, since~$\|u_\star\|_{L^{p+1}(\R^n_+)}=\|u\|_{L^{p+1}(\R^n_+)}=1$, minimality entails that~$\|u_\star\|_{s,+}=\|u\|_{s,+}$.
Thus, the difference in~\eqref{rearrangement:regional} vanishes for almost every~$(x_n,y_n)\in(0,+\infty)^2$.

Now, let~$y_n>0$ be such that~$u(\cdot,y_n)\in L^2(\R^{n-1})$ and the equality in~\eqref{rearrangement:regional} holds for almost every~$x_n>0$. The equality case in Riesz's rearrangement inequality
in~\cite[Theorem~3.9]{MR1817225}, applied to the strictly decreasing kernel~$(|x'-y'|^2+|x_n-y_n|^2)^{-(n+2s)/2}$ for~$x_n\ne y_n$, shows that~$u(\cdot,x_n)$ and~$u(\cdot,y_n)$ are translates of their rearrangements by the same tangential vector, say~$\varpi$,
and the nonzero slice~$u(\cdot,y_n)$ fixes the vector~$\varpi$ for almost every~$x_n$. 

Thus, translating~$u$ by~$\varpi$ and using the interior continuity ensured by~\eqref{OFCLA}, we get that~$u(\cdot,x_n)=u_\star(\cdot,x_n)$ for every~$x_n>0$. Moreover, the Neumann formula extends the rotational symmetry to~$x_n<0$, and continuity from~$\R^n_+$ extends the symmetry to~$x_n=0$. The proof of~\eqref{REARR:in} is thereby complete.}

\noindent{\bf Step 4. Radial monotonicity and proof of~\eqref{eqRedde}.}
We know from Step~3 and~\eqref{OFCLA} that, for all~$\rho>0$ and~$x_n\in\R$,
\begin{equation}\label{ANC:nos} \partial_\rho\overline u(\rho,x_n)\le0.\end{equation}

Given~$i\in\{1,\ldots,n-1\}$,
let~$H_i:=\{x_i>0\}$. 
Recalling the notation introduced in~\eqref{defAr7436983480},
the Neumann condition reads, for~$x_n<0$, 
\begin{equation*}
u(x)=\frac1{A(x)}\int_{\R^n_+}\frac{u(y)}{|x-y|^{n+2s}}\,dy.
\end{equation*}
This and~\eqref{OFCLA}
allow us to define~$w:=-\partial_{x_i}u$ in~$\R^n\setminus\{x_n=0\}$.

Thus, for all~$x\in H_i$ such that~$ x_n\ne0$,
the radial symmetry in~\eqref{REARR:in}
and the monotonicity in~\eqref{ANC:nos} lead to
\begin{equation}\label{cbjewgyrwtr87654}
w(x-2x_i e_i)=-w(x)\qquad{\mbox{and}}\qquad
w(x)=-\frac{x_i}{|x'|}\partial_\rho\overline u(|x'|,x_n)\ge0
.\end{equation}
Hence, for almost every~$(x_1,\dots,x_{i-1},x_{i+1},\dots,x_n)$, 
\begin{equation*}
\int_{\R}|w(x)|\,{dx_i}
=2\left(u(x_1,\dots,x_{i-1},0,x_{i+1},\dots,x_n)-\lim_{x_i\to+\infty}u(x)\right)\le2\|u\|_{L^\infty(\R^n)}.
\end{equation*}
Consequently,
\begin{equation}\label{TH:by1-1}\begin{split}&
\int_{\R^n}\frac{|w(y)|}{(1+|y|)^{n+2s}}\,dy
\le
\int_{\R^n}\frac{|w(y)|}{(1+|(y_1,\dots,y_{i-1},y_{i+1},\dots,y_n)|)^{n+2s}}\,dy\\&\qquad=
\int_{\R^{n-1}}\left(\int_\R |w(y)|\,dy_i\right)
\frac{dy_1\dots dy_{i-1}\,dy_{i+1}\dots dy_n}{(1+|(y_1,\dots,y_{i-1},y_{i+1},\dots,y_n)|)^{n+2s}}
\\&\qquad\le
2\|u\|_{L^\infty(\R^n)}\int_{\R^{n-1}}\frac{dz}{(1+|z|)^{n+2s}}<+\infty.
\end{split}
\end{equation}

On this account, we are allowed to perform tangential integration by parts and, since~$A$ is independent of~$x_i$, we obtain
from~\eqref{rearrangement:neumann-average} that
\begin{equation}\label{TH:by1-2}
\begin{split}
w(x)
&=-\frac1{A(x)}\int_{\R^n_+}u(y)\partial_{x_i}|x-y|^{-n-2s}
\,dy\\
&=\frac1{A(x)}\int_{\R^n_+}u(y)\partial_{y_i}{|x-y|^{-n-2s}}\,dy\\&
=\frac1{A(x)}\int_{\R^n_+}\frac{w(y)}{|x-y|^{n+2s}}\,dy.\end{split}
\end{equation}
As a result, for all~$x\in\R^n_-$,
\begin{equation*}
\mathcal N_s w(x)=c_{n,s}\left(A(x)w(x)-\int_{\R^n_+}\frac{w(y)}{|x-y|^{n+2s}}\,dy\right)=0.
\end{equation*}

Now, to differentiate the equation
in a rigorous way, let~$\varphi\in C^\infty_0(\R^n_+)$. We deduce from~\eqref{TH:by1-1} and~\eqref{TH:by1-2} that
\begin{align*}&
\int_{\R^n}w(-\Delta)^s\varphi\,dx
=\int_{\R^n}u\,\partial_{x_i}(-\Delta)^s\varphi\,dx
=\int_{\R^n}u(-\Delta)^s(\partial_{x_i}\varphi)\,dx\\
&\qquad=\int_{\R^n_+}(\nu_su^p-u)\partial_{x_i}\varphi\,dx
=\int_{\R^n_+}(\nu_spu^{p-1}-1)w\varphi\,dx
\end{align*}
and therefore~$
(-\Delta)^sw+w=\nu_spu^{p-1}w$ pointwise in~$\R^n_+$.

Suppose now, for the sake of contradiction, that~$w(x_0)=0$ at some~$x_0\in H_i\cap\R^n_+$. Then, by the nonnegativity of~$w$ (recall~\eqref{cbjewgyrwtr87654}), it follows that~$\nabla w(x_0)=0$, and therefore~$w(y)=O(|y-x_0|^2)$ near~$x_0$. 

On this account, the integrals in the following computations
converge absolutely. More precisely, we can split~$\R^n$ into~$H_i$ and
its complement and, using the antisymmetry of~$w$ in~\eqref{cbjewgyrwtr87654}, we find that
\begin{equation}\label{DZcqQ}
\begin{split}
0&=(-\Delta)^sw(x_0)\\
&=-c_{n,s}\int_{H_i}\left(\frac{w(y)}{|x_0-y|^{n+2s}}+\frac{w(y-2y_i e_i)}{|x_0-y+2y_i e_i|^{n+2s}}
\right)\,dy\\
&=-c_{n,s}\int_{H_i}w(y)\left(\frac{1}{|x_0-y|^{n+2s}}-\frac{1}{|x_0-y+2y_i e_i|^{n+2s}}\right)\,dy.
\end{split}\end{equation}
Now, for every~$y\in H_i$, we have that
\begin{equation*}
\begin{split}
&{|x_0-y+2y_i e_i|^2}-|x_0-y|^2=4x_{0,i} y_i>0\\
\mbox{and}\qquad&\frac{1}{|x_0-y|^{n+2s}}>\frac{1}{|x_0-y+2y_i e_i|^{n+2s}}.
\end{split}
\end{equation*}
Therefore, since~$w\ge0$ in~$H_i$, we deduce from~\eqref{DZcqQ}
that~$w=0$ a.e. in~$H_i$. The antisymmetry property in~\eqref{cbjewgyrwtr87654} then gives that~$\partial_{x_i}u\equiv0$ in~$\R^n_+$. This observation, together with~\eqref{REPOSIt},
contradicts the fact that~$u\in L^{p+1}(\R^n_+)$.

Therefore~$w>0$ in~$H_i\cap\R^n_+$. Taking~$x=\rho e_i+t e_n$, with~$\rho$, $t>0$, gives
\begin{equation*}
-\partial_\rho\overline u(\rho,t)=w(\rho e_i+t e_n)>0,
\end{equation*}
which proves~\eqref{eqRedde}.
\end{proof}

\section{Localization of the maximum of the ground state}\label{LOC:MAXST1}
Here we show that the location of the maximum of the ground state solution remains bounded when~$s\nearrow1$. More specifically, if~$u_s$ is the normalized ground state solution of problem~\eqref{problem_nu} constructed in Theorems~\ref{EXIST:GRS} and~\ref{ADDITIO:PR},
we have that~$u_s(0,t)\to0$ as $t\to+\infty$,
{by~\eqref{OFCLA} and the~$L^{p+1}$ normalization in~\eqref{MIN:GRO:ST}}.
We also know from the nonlocal Neumann condition that~$u_s(x)\le\|u_s\|_{L^\infty(\R^n_+)}$ for all~$x\in\R^n$ and,
in light of~\eqref{REARR:in}, \eqref{OFCLA}, and~\eqref{eqRedde}, we have that
\begin{equation}\label{MAXDE}
{M}_s:=\|u_s\|_{L^\infty(\R^n)}=
\|u_s\|_{L^\infty(\R^{n-1}\times[0,+\infty))}=u_s({x}_s),\end{equation}
with~${x}_s=(0,\dots,0,{m}_s)$ for some~${m}_s\in[0,+\infty)$.

We now prove the following result:

\begin{proposition}\label{DOVEMI}
We have that
$$ \limsup_{s\nearrow1} m_{s}<+\infty.$$
\end{proposition}

The proof of Proposition~\ref{DOVEMI} relies on the following uniform bounds:

\begin{lemma}\label{SI-en-b} There exists a positive constant~$C^\star_{n,p}$ such that~$\nu_s\le
C^\star_{n,p}$.
\end{lemma}

\begin{proof} Let~$\phi\in C^\infty_c((1,2)^n, [0,+\infty))$
with~$\|\phi\|_{L^{p+1}(\R^n_+)}=1$. Thus,
\begin{eqnarray*}
\nu_s-\|\phi\|^2_{L^{2}(\R^n_+)}&\le&
\frac{c_{n,s}}2\iint_{Q(\R^n_+)}\frac{(\phi(x)-\phi(y))^2}{|x-y|^{n+2s}}\,dx\,dy\\&
\le& c_{n,s}\iint_{(1,2)^n\times\R^{n}}\frac{(\phi(x)-\phi(y))^2}{|x-y|^{n+2s}}\,dx\,dy\\&
\le& 4c_{n,s}\,\|\phi\|_{C^1(\R^n)}^2\iint_{(1,2)^n\times\R^{n}}\frac{\min\{1,|x-y|^2\}}{|x-y|^{n+2s}}\,dx\,dy\\&\le& C_n c_{n,s}\,\|\phi\|_{C^1(\R^n)}^2\left(\frac1{1-s}+\frac1s\right),
\end{eqnarray*}
where~$Q(\R^n_+)=\R^{2n}\setminus(\R^n_-\times\R^n_-)$.
Therefore, the desired result follows from the asymptotics of~$c_{n,s}$ (see e.g.~\cite[Section~4]{MR2944369}).
\end{proof}

\begin{lemma}\label{3e4rfOSJL:LEM}
Let~$u_s$ be a positive normalized ground state of~\eqref{problem_nu}, as in Theorems~\ref{EXIST:GRS} and~\ref{ADDITIO:PR}. 

Then there exist~$s_0\in(0,1)$ and~$C>0$, depending only on~$n$, $p$ and~$s_0$, such that
$$\sup_{s\in[s_0,1)} \|u_s\|_{L^\infty(\R^n_+)}\le C.$$
\end{lemma}

\begin{proof}
We take~$s_0\in(0,1)$ with~$p+1<2_{s_0}^*$ and set~$\kappa:=2_{s_0}^*/(p+1)>1$. 
The proof is now a revisitation of Step~2 in the proof of Theorem~\ref{ADDITIO:PR}, checking that the constant in~\eqref{ground:moser-bound} is independent of~$s$ for any~$s\ge s_0$.

For this, we observe that~$2<2_{s_0}^*\le2_s^*$ for every~$s\ge s_0$, so we can apply~\eqref{ground:sobolev} with the constant~$\widetilde{C}$ depending only on~$n$. 

Moreover, recalling the notation in~\eqref{bcjkegity4gw78yt348},
we see that~$A=\max\{1,\widetilde{C}\nu_s\}$ is bounded independently of~$s$, thank to Lemma~\ref{SI-en-b}. Hence, one can use~\eqref{ground:moser-bound} and obtain the desired uniform bound.
\end{proof}

It is also useful to refine the ground state energy comparison in~\eqref{EWNS:02oe-034.P}
to obtain a uniform bound in~$s$. Specifically, the ground state energy in the halfspace
{is strictly smaller than the full-space ground state energy, uniformly} as~$s\nearrow1$, as pointed out in the next result:

\begin{lemma}\label{SJMSD:ALmdloi} Let~$\nu_s^\star$ be as in~\eqref{def_nustar}.
Then, there exists~$c_{n,p}>0$ such that
\begin{equation*}
\lim_{s\nearrow1}\nu_s^\star
\ge2^{1-\frac2{p+1}}\limsup_{s\nearrow1}\nu_s\qquad{\mbox{and}}\qquad
\liminf_{s\nearrow1}(\nu_s^\star-\nu_s)\ge c_{n,p}.
\end{equation*}
\end{lemma}

\begin{proof}
We have that
\begin{equation}\label{halfspace:fixed-test}
\begin{split}\lim_{s\nearrow1}
\|\Phi\|_{s,+}^2
&=\lim_{s\nearrow1}\|\Phi\|_{H^s(\R^n)}^2
-\frac{c_{n,s}}2\iint_{\R^n_-\times\R^n_-}
\frac{|\Phi(x)-\Phi(y)|^2}{|x-y|^{n+2s}}\,dx\,dy
-\int_{\R^n_-}\Phi^2\,dx\\
&=
\|\Phi\|_{H^1(\R^n)}^2
-\int_{\R^n_-}(|\nabla\Phi|^2+\Phi^2)\,dx\\&
=\|\Phi\|_{H^1(\R^n_+)}^2.
\end{split}
\end{equation}
In~\eqref{halfspace:fixed-test}, the limit in~$\R^n$ follows by Plancherel's theorem and the dominated convergence theorem, while the limit in~$\R^n_-$ follows from~\cite[Theorem~1]{MR4093822} and~\cite[Corollary~4.2(i)]{MR2944369}
(of course, $\R^n_-$ is an extension domain).

{Let
\begin{equation}\label{nu1:whole-space}
\nu_1^\star:=\inf_{u\in H^1(\R^n)\setminus\{0\}}
\frac{\|u\|_{H^1(\R^n)}^2}{\|u\|_{L^{p+1}(\R^n)}^2}.
\end{equation}
Let~$\omega_1$ be
a positive radial minimizer, normalized by~$\|\omega_1\|_{L^{p+1}(\R^n)}=1$, as in~\cite{MR695535}. Then, by symmetry,
\begin{equation*}
\|\omega_1\|_{L^{p+1}(\R^n_+)}^2=2^{-\frac{2}{p+1}}
\qquad{\mbox{and}}\qquad
\|\omega_1\|_{H^1(\R^n_+)}^2=\frac12\nu_1^\star.
\end{equation*}
Therefore, the variational definition~\eqref{MIN:GRO:ST} and~\eqref{halfspace:fixed-test} imply that
\begin{equation*}
\limsup_{s\nearrow1} 2^{-\frac{2}{p+1}}\nu_s
\le \limsup_{s\nearrow1}\|\omega_1\|_{s,+}^2
=\frac12\nu_1^\star.
\end{equation*}
Since~$\nu_s^\star\to\nu_1^\star$ by~\cite[Corollary~2.5]{Fall-Valdinoci}, we have that
\begin{equation*}
\lim_{s\nearrow1}\nu_s^\star=\nu_1^\star
\ge2^{1-\frac{2}{p+1}}\limsup_{s\nearrow1}\nu_s.
\end{equation*}
Consequently, recalling also~\eqref{NONZGSI},
\begin{equation*}
\liminf_{s\nearrow1}(\nu_s^\star-\nu_s)
=\nu_1^\star-\limsup_{s\nearrow1}\nu_s
\ge\left(2^{1-\frac{2}{p+1}}-1\right)c^\star_{n,p}>0.
\end{equation*}
The proof is thereby complete.}
\end{proof}

With this, we can now proceed with the proof of Proposition~\ref{DOVEMI}.

\begin{proof}[Proof of Proposition~\ref{DOVEMI}] Suppose not. Then, there {exists} a sequence~$s_k\nearrow1$ such that
\begin{equation}\label{AMoawdSJhyA}
{\mbox{$m_{s_k}\to+\infty$ as~$k\to+\infty$.}}\end{equation}
To track the energy and mass near the maximum, we define
\begin{equation*}\begin{aligned}
P_k(R)&:=\frac{c_{n,s_k}}2\iint_{{\R^n\times(\R^n_+\cap B_R(x_{s_k}))}}\frac{(u_{s_k}(x)-u_{s_k}(y))^2}{|x-y|^{n+2{s_k}}}\,dx\,dy\\
&\qquad+\int_{\R^n_+\cap B_R(x_{s_k})} \big(u_{s_k}^2(x)+u_{s_k}^{p+1}(x)\big)\,dx.
\end{aligned}\end{equation*}
We {observe} that~$P_k$ is nondecreasing in~$R$ and~$P_k(R)\le\nu_{s_k}+1$, which is bounded uniformly in~$k$, thanks to Lemma~\ref{SI-en-b}.
Thus, up to a subsequence, we can apply Helly's selection theorem and find a nondecreasing function~$P_\infty$ such that, for almost all~$R>0$,
$$ P_\infty(R)=\lim_{k\to+\infty}P_k(R).$$

Let now~$m\in\N\setminus\{0\}$. Since
$$ \lim_{R\to+\infty} P_\infty(R)=\sup_{(0,+\infty)}P_\infty<+\infty,$$
we can find~$R_m>0$ sufficiently large such that, for all~$\rho$, $R\ge R_m$, we have that\begin{equation*}|P_\infty(R)-P_\infty(\rho)|\le\frac1m.\end{equation*}
Without loss of generality, we can also suppose that~$R_{m+1}>R_m+2m$. In this way, we find~$\rho_m\in[R_m,R_m+1]$ such that
$$\frac1m\ge|P_\infty(\rho_{m+1})-P_\infty(\rho_m)|
=\lim_{k\to+\infty}|P_k(\rho_{m+1})-P_k(\rho_m)|.$$

By this and~\eqref{AMoawdSJhyA}, we can find~$k_m\in\N$ such that for all~$k\ge k_m$ we have that~$m_{s_k}\ge \rho_m+2m$ and
\begin{equation}\label{1r2e3d7OJASLSsIKS}\begin{split}\frac2m&\ge|P_k(\rho_{m+1})-P_k(\rho_m)|
\\&=\frac{c_{n,s_k}}2\iint_{{\R^n\times\big(\R^n_+\cap(B_{\rho_{m+1}}(x_{s_k})\setminus B_{\rho_m}(x_{s_k}))\big)}}\frac{(u_{s_k}(x)-u_{s_k}(y))^2}{|x-y|^{n+2{s_k}}}\,dx\,dy\\&\qquad\qquad
+\int_{\R^n_+\cap(B_{\rho_{m+1}}({x_{s_k}})\setminus B_{\rho_{m}} ({x_{s_k}}))} \big(u_{s_k}^2(x)+u_{s_k}^{p+1}(x)\big)\,dx.
\end{split}\end{equation}

We now consider the annulus
$$A_m:=B_{\rho_m+m}{(x_{s_{k_m}})}\setminus B_{\rho_m}{(x_{s_{k_m}})}.$$ Since~$m_{s_{k_m}}\ge\rho_m+2m$, we have that~$B_{\rho_m+m}{(x_{s_{k_m}})}\subset\R^n_+$.
Let also~$\tau_m\in C^\infty(\R^n,[0,1])$ with~$\tau_m=1$ in~$B_{\rho_m}(x_{s_{k_m}})$ and~$\tau_m=0$ in~$\R^n\setminus B_{\rho_m+m}{(x_{s_{k_m}})}$, such that~$|\nabla \tau_m|\le\frac{C}m$. Let~$v_m:=u_{s_{k_m}}\tau_m$ and~$w_m:=u_{s_{k_m}}-v_m=(1-\tau_m)u_{s_{k_m}}$.

Up to a subsequence in~$m$, we may assume that~$\nu_{s_{k_m}}$ and the bounded sequences
$$\int_{\R^n_+}v_m^q\,dx\qquad{\mbox{and}}\qquad\int_{\R^n_+}w_m^q\,dx,$$ for~$q\in\{2,p+1\}$, converge.

Additionally, for~$q\in\{2,p+1\}$,
we have that~$ u_{s_{k_m}}^q-v_m^q-w_m^q=0$ in~$B_{\rho_m}{(x_{s_{k_m}})}$ and outside~$ B_{\rho_m+m}{(x_{s_{k_m}})}$.
Moreover, 
$$\rho_{m+1}\ge R_{m+1}\ge R_m+2m\ge \rho_m+2m-1\ge \rho_m+m.$$
Therefore, by~\eqref{1r2e3d7OJASLSsIKS},
\begin{equation*}\begin{split}&
\int_{\R^n_+} |u_{s_{k_m}}^q-v_m^q-w_m^q|\,dx=
\int_{ \R^n_+\cap (B_{\rho_m+m}{(x_{s_{k_m}})}\setminus B_{\rho_m}{(x_{s_{k_m}})})} |u_{s_{k_m}}^q-v_m^q-w_m^q|\,dx\\&\quad\le
C \int_{ \R^n_+\cap (B_{\rho_{m+1}}{(x_{s_{k_m}})}\setminus B_{\rho_m}{(x_{s_{k_m}})})} u_{s_{k_m}}^q\,dx\le\frac{2C}m,
\end{split}\end{equation*}
where~$C$ is a positive constant possibly varying from line to line,
and accordingly
\begin{equation}\label{DE:st0}
\lim_{m\to+\infty}\int_{\R^n_+} u_{s_{k_m}}^q(x)\,dx=\lim_{m\to+\infty}\int_{\R^n_+}v_m^q(x)\,dx+
\int_{\R^n_+}w_m^q(x)\,dx.
\end{equation}

Next, we use the notation~$Q(\R^n_+)=\R^{2n}\setminus(\R^n_-\times\R^n_-)$ and
we claim that
\begin{equation}\label{DE:st1}
\lim_{m\to+\infty}\frac{c_{n,s_{k_m}}}2
\iint_{Q(\R^n_+)}\frac{\big(v_m(x)-v_m(y)\big)\big(w_m(x)-w_m(y)\big)}{|x-y|^{n+2{s_{k_m}}}}\,dx\,dy
=0.\end{equation}
To check this, we point out that
\begin{eqnarray*}
v_m(x)-v_m(y)&=&u_{s_{k_m}}(x)\tau_m(x)-u_{s_{k_m}}(y)\tau_m(y)\\&=&
u_{s_{k_m}}(x)\big(\tau_m(x)-\tau_m(y)\big)
+\tau_m(y)\big(u_{s_{k_m}}(x)-u_{s_{k_m}}(y)\big)
\end{eqnarray*}
and similarly, replacing~$\tau_m$ with~$1-\tau_m$,
\begin{eqnarray*}
&&w_m(x)-w_m(y)=u_{s_{k_m}}(x)\big(\tau_m(y)-\tau_m(x)\big)
+\big(1-\tau_m(y)\big)\big(u_{s_{k_m}}(x)-u_{s_{k_m}}(y)\big).
\end{eqnarray*}
As a result, given~$\delta\in{(0,\tfrac12)}$,
\begin{eqnarray*}&&
\big|v_m(x)-v_m(y)\big|\,\big|w_m(x)-w_m(y)\big|\\&&\quad\le
u_{s_{k_m}}^2(x)\big(\tau_m(x)-\tau_m(y)\big)^2\\&&\qquad+
\big(1-\tau_m(y)\big)\big|\tau_m(x)-\tau_m(y)\big|\,|u_{s_{k_m}}(x)|\,
\big|u_{s_{k_m}}(x)-u_{s_{k_m}}(y)\big|\\&&\qquad\qquad
+\tau_m(y) \big|\tau_m(y)-\tau_m(x)\big|\,
|u_{s_{k_m}}(x)|\,
\big|u_{s_{k_m}}(x)-u_{s_{k_m}}(y)\big|\\&&\qquad\qquad
+\tau_m(y) \big(1-\tau_m(y)\big)
\big(u_{s_{k_m}}(x)-u_{s_{k_m}}(y)\big)^2\\&&\quad=
u_{s_{k_m}}^2(x)\big(\tau_m(x)-\tau_m(y)\big)^2+
\big|\tau_m(x)-\tau_m(y)\big|\,|u_{s_{k_m}}(x)|\,
\big|u_{s_{k_m}}(x)-u_{s_{k_m}}(y)\big|\\&&\qquad\qquad
+\tau_m(y) \big(1-\tau_m(y)\big)
\big(u_{s_{k_m}}(x)-u_{s_{k_m}}(y)\big)^2\\&&\quad\le\frac{C}\delta
u_{s_{k_m}}^2(x)\big(\tau_m(x)-\tau_m(y)\big)^2
+C\delta\big(u_{s_{k_m}}(x)-u_{s_{k_m}}(y)\big)^2\\&&\qquad\qquad
+\tau_m(y) \big(1-\tau_m(y)\big)
\big(u_{s_{k_m}}(x)-u_{s_{k_m}}(y)\big)^2.
\end{eqnarray*}

Since the integrand in~\eqref{DE:st1} vanishes
when both~$x$ and~$y$ lie outside~$B_{\rho_m+m}{(x_{s_{k_m}})}$, we obtain
\begin{equation}\label{gTGBN8ba20f93-1}
\begin{split}&\frac{c_{n,s_{k_m}}}2
\iint_{Q(\R^n_+)}\frac{\big|v_m(x)-v_m(y)\big|\,\big|w_m(x)-w_m(y)\big|}{|x-y|^{n+2{s_{k_m}}}}\,dx\,dy\\&\quad\le{c_{n,s_{k_m}}}
\iint_{B_{\rho_m+m}{(x_{s_{k_m}})}\times\R^n}\frac{\big|v_m(x)-v_m(y)\big|\,\big|w_m(x)-w_m(y)\big|}{|x-y|^{n+2{s_{k_m}}}}\,dx\,dy\\&\quad\le
\frac{c_{n,s_{k_m}}}{\delta}
\iint_{B_{\rho_m+m}{(x_{s_{k_m}})}\times\R^n}\frac{u_{s_{k_m}}^2(x)\big(\tau_m(x)-\tau_m(y)\big)^2}{|x-y|^{n+2{s_{k_m}}}}\,dx\,dy
+2\delta\nu_{s_{k_m}}\\&\qquad\qquad+ c_{n,s_{k_m}}
\iint_{B_{\rho_m+m}{(x_{s_{k_m}})}\times\R^n}\frac{\tau_m(y) \big(1-\tau_m(y)\big)
\big(u_{s_{k_m}}(x)-u_{s_{k_m}}(y)\big)^2}{|x-y|^{n+2{s_{k_m}}}}\,dx\,dy.\end{split}
\end{equation}

We now observe that
\begin{equation}\label{gTGBN8ba20f93-2}\begin{split}&
\iint_{B_{\rho_m+m}{(x_{s_{k_m}})}\times\R^n}\frac{u_{s_{k_m}}^2(x)\big(\tau_m(x)-\tau_m(y)\big)^2}{|x-y|^{n+2{s_{k_m}}}}\,dx\,dy\\&\quad
\le
C\iint_{B_{\rho_m+m}{(x_{s_{k_m}})}\times\R^n} u_{s_{k_m}}^2(x)
\min\left\{1,\frac{|x-y|^2}{m^2}\right\}\frac{dx\,dy}{|x-y|^{n+2{s_{k_m}}}}
\\&\quad
\le \frac{C}{m^{2s_{k_m}}}\left(\frac1{1-s_{k_m}}+\frac1{s_{k_m}}\right){\int_{B_{\rho_m+m}{(x_{s_{k_m}})}}} u_{s_{k_m}}^2(x)\,dx\\&\quad
\le \frac{C}{c_{n,s_{k_m}}\,m^{2s_{k_m}}},
\end{split}
\end{equation}up to relabeling~$C$
(the last step being a consequence of the asymptotics of~$c_{n,s}$, see e.g.~\cite[Section~4]{MR2944369}).

Also, owing to~\eqref{1r2e3d7OJASLSsIKS},
\begin{eqnarray*}&&
c_{n,s_{k_m}}
\iint_{B_{\rho_m+m}{(x_{s_{k_m}})}\times\R^n}\frac{\tau_m(y) \big(1-\tau_m(y)\big)
\big(u_{s_{k_m}}(x)-u_{s_{k_m}}(y)\big)^2}{|x-y|^{n+2{s_{k_m}}}}\,dx\,dy\\&&\qquad\le
c_{n,s_{k_m}}
\iint_{{\R^n\times A_m}}\frac{
\big(u_{s_{k_m}}(x)-u_{s_{k_m}}(y)\big)^2}{|x-y|^{n+2{s_{k_m}}}}\,dx\,dy\le\frac4m 
,\end{eqnarray*}
which is infinitesimal as~$m\to+\infty$.

Therefore, inserting this information and~\eqref{gTGBN8ba20f93-2} into~\eqref{gTGBN8ba20f93-1}, and recalling Lemma~\ref{SI-en-b}, we find that
\begin{equation*}
{\limsup_{m\to+\infty}}\,\frac{c_{n,s_{k_m}}}2
\iint_{Q(\R^n_+)}\frac{\big|v_m(x)-v_m(y)\big|\,\big|w_m(x)-w_m(y)\big|}{|x-y|^{n+2{s_{k_m}}}}\,dx\,dy\le C\delta.\end{equation*}
This completes the proof of~\eqref{DE:st1}, up to taking~$\delta$ as small as we wish.

{F}rom~\eqref{DE:st1} we arrive at
\begin{equation}\label{DE:st2}\begin{split}
&\!\!\!\!\!\lim_{m\to+\infty}
\frac{c_{n,s_{k_m}}}2\iint_{Q(\R^n_+)}\frac{(u_{s_{k_m}}(x)-u_{s_{k_m}}(y))^2}{|x-y|^{n+2{s_{k_m}}}}\,dx\,dy\\&
=\lim_{m\to+\infty}
\frac{c_{n,s_{k_m}}}2\iint_{Q(\R^n_+)}\frac{\big((v_m+w_m)(x)-(v_m+w_m)(y)\big)^2}{|x-y|^{n+2{s_{k_m}}}}\,dx\,dy
\\&=\lim_{m\to+\infty}
\bigg({\frac{c_{n,s_{k_m}}}2\iint_{Q(\R^n_+)}\frac{(v_m(x)-v_m(y))^2}{|x-y|^{n+2{s_{k_m}}}}\,dx\,dy}\\&\quad\qquad\qquad
+\frac{c_{n,s_{k_m}}}2\iint_{Q(\R^n_+)}\frac{(w_m(x)-w_m(y))^2}{|x-y|^{n+2{s_{k_m}}}}\,dx\,dy\bigg).
\end{split}\end{equation}

Therefore, using~\eqref{DE:st0} and~\eqref{DE:st2} together with the convexity inequality in~\eqref{WEHi5oFI}, we conclude that
\begin{eqnarray*}
\lim_{m\to+\infty} {\nu_{s_{k_m}}}&=&\lim_{m\to+\infty}
\bigg(\frac{c_{n,s_{k_m}}}2\iint_{Q(\R^n_+)}\frac{(u_{s_{k_m}}(x)-u_{s_{k_m}}(y))^2}{|x-y|^{n+2{s_{k_m}}}}\,dx\,dy
+\int_{\R^n_+} u_{s_{k_m}}^2(x)\,dx\bigg)
\\&=&\lim_{m\to+\infty}
\bigg(\frac{c_{n,s_{k_m}}}2\iint_{Q(\R^n_+)}\frac{(v_m(x)-v_m(y))^2}{|x-y|^{n+2{s_{k_m}}}}\,dx\,dy
+\int_{\R^n_+} v_m^2(x)\,dx\\&&\quad
+\frac{c_{n,s_{k_m}}}2\iint_{Q(\R^n_+)}\frac{(w_m(x)-w_m(y))^2}{|x-y|^{n+2{s_{k_m}}}}\,dx\,dy
+\int_{\R^n_+} w_m^2(x)\,dx\bigg)\\&\ge&
\lim_{m\to+\infty} \nu_{s_{k_m}}^\star \|v_m\|^2_{L^{p+1}(\R^n_+)}+
{\nu_{s_{k_m}}}\|w_m\|^2_{L^{p+1}(\R^n_+)}\\&=&
\lim_{m\to+\infty} {\nu_{s_{k_m}}} \Big(\|v_m\|^2_{L^{p+1}(\R^n_+)}+\|{w_m}\|^2_{L^{p+1}(\R^n_+)}\Big)+
(\nu_{s_{k_m}}^\star-{\nu_{s_{k_m}}})\|v_m\|^2_{L^{p+1}(\R^n_+)}
\\&\ge&
\lim_{m\to+\infty} {\nu_{s_{k_m}}}
\Big( \|v_m\|^{p+1}_{L^{p+1}(\R^n_+)}+\|w_m\|^{p+1}_{L^{p+1}(\R^n_+)}\Big)^{\frac2{p+1}}+(\nu_{s_{k_m}}^\star-{\nu_{s_{k_m}}})\|v_m\|^2_{L^{p+1}(\R^n_+)}\\&=&
\lim_{m\to+\infty} {\nu_{s_{k_m}}}\big(\|u_{s_{k_m}}\|^{p+1}_{L^{p+1}(\R^n_+)}\big)^{\frac2{p+1}}+(\nu_{s_{k_m}}^\star-{\nu_{s_{k_m}}})\|v_m\|^2_{L^{p+1}(\R^n_+)}\\&=&\lim_{m\to+\infty} {\nu_{s_{k_m}}}+(\nu_{s_{k_m}}^\star-{\nu_{s_{k_m}}})\|v_m\|^2_{L^{p+1}(\R^n_+)}
\end{eqnarray*}
and therefore
$$\lim_{m\to+\infty}(\nu_{s_{k_m}}^\star-{\nu_{s_{k_m}}})\|v_m\|^2_{L^{p+1}(\R^n_+)}=0.$$

This, together with Lemma~\ref{SJMSD:ALmdloi}, gives that
\begin{equation}\label{OJSLMDSKPDMON2e034tujh4}
\lim_{m\to+\infty}\|v_m\|_{L^{p+1}(\R^n_+)}=0.
\end{equation}

In the rest of the proof we take~$s=s_{k_m}$ and~$m$ sufficiently large such that~$B_2(x_s)\subset\R^n_+$ and~$\rho_m\ge1$.
It is also interesting to observe that, at the maximal point~$x_s$, we have that~$(-\Delta)^s u_s(x_s)\ge0$ and accordingly, recalling~\eqref{MAXDE},
$$ M_s=u_s(x_s)\le(-\Delta)^s u_s(x_s)+u_s(x_s)=
\nu_s u_s^p(x_s)=\nu_s M_s^p.$$
Therefore, exploiting Lemma~\ref{SI-en-b},
\begin{equation}\label{LOSN} M_s\ge \left(\frac1{\nu_s}\right)^{\frac1{p-1}}\ge\frac1{C_0},\end{equation}
for some universal constant~$C_0$, depending only on~$n$ and~$p$.

Now, by Lemma~\ref{3e4rfOSJL:LEM} and the nonlocal Neumann condition, we deduce that~$\|u_s\|_{L^\infty(\R^n)}$ is bounded uniformly in~$s$. Thus, by
interior fractional elliptic regularity, it follows that~${\|u_s\|_{C^{\alpha}(B_1(x_s))}}$ is bounded uniformly in~$s$, for some~$\alpha\in(0,1)$ (see e.g.~\cite[Theorem~12.1]{MR2494809}).
As a result, recalling~\eqref{LOSN}, there exists~$r_0\in(0,1)$, independent {of}~$s$, such that, for all~$x\in B_{r_0}(x_s)$,
$$ u_s(x)\ge u_s(x_s)- |u_s(x)-u_s(x_s)|=M_s- |u_s(x)-u_s(x_s)|\ge\frac1{C_0}-
|u_s(x)-u_s(x_s)|\ge\frac1{2C_0}.
$$
Consequently,
$$ \|v_m\|_{L^{p+1}(\R^n_+)}^{p+1}\ge \int_{B_{r_0}(x_s)} v_m^{p+1}(x)\,dx=\int_{B_{r_0}(x_s)} u_{s_{k_m}}^{p+1}(x)\,dx
\ge\frac{|B_{r_0}|}{{(2C_0)^{p+1}}}
,$$
but this is in contradiction with~\eqref{OJSLMDSKPDMON2e034tujh4}.\end{proof}

\section{Asymptotic analysis}\label{asymptotics}
Let
$$\mathcal{M}_s:=\{u_s\in H^s_+ \mbox{ ground state solutions of } \eqref{problem_nu} \mbox{ for }\nu_s\}.$$

Let us denote
\begin{equation}\label{nu1}
\nu_1:=\inf_{{u\in H^1(\R^n_+)}\atop{\|u\|_{L^{p+1}(\R^n_+)}=1}}\|u\|_{H^1(\R^n_+)}^2,
\end{equation}
and, recalling the definition of~$\nu_1^\star$ in~\eqref{nu1:whole-space},
\begin{equation*}\begin{split}
\mathcal{M}_1&:=\{u\in H^1(\R^n_+)\;\mbox{ such that }\;\|u\|_{L^{p+1}(\R^n_+)}=1\;\mbox{ and }\;\|u\|^2_{H^1(\R^n_+)}=\nu_1\}\\{\mbox{and }}\quad
\mathcal{M}_1^*&:=\{u\in H^1(\R^n)\;\mbox{ such that }\;\|u\|_{L^{p+1}(\R^n)}=1\;\mbox{ and }\;\|u\|^2_{H^1(\R^n)}=\nu_1^\star\}.
\end{split}\end{equation*}

\begin{remark}\label{relationnu}
If $u\in \mathcal{M}_1^*$ with maximum value at any point of the form $(x',0)$, then $2^{\frac{1}{p+1}}u|_{\R^n_+}\in \mathcal{M}_1$. Likewise, if $u\in \mathcal{M}_1$ and $\tilde{u}$ denotes its even extension to the whole space, then $2^{\frac{-1}{p+1}}\tilde{u}$ belongs to $\mathcal{M}_1^*$. As a {consequence}
$$\nu_1=2^{\frac{2}{p+1}-1}\nu_1^\star.$$
\end{remark}

We first prove a lower-semicontinuity estimate.

\begin{proposition}\label{prop:liminf}
Let $\Omega\subset\R^n$ be an open set and $u_s\in L^2(\Omega)$ be such that, as~$s\nearrow1$,
$$u_s\to u\quad\text{in }L^2_{\rm loc}(\Omega)\qquad\text{for some }u\in H^1(\Omega).$$
Then
\begin{equation*}\label{eq:liminf}
\liminf_{s\nearrow1}\frac{c_{n,s}}{2}\iint_{\Omega\times\Omega}\frac{(u_s(x)-u_s(y))^2}{|x-y|^{n+2s}}\,dx\,dy\ \ge\ \int_\Omega|\nabla u|^2\,dx.
\end{equation*}
\end{proposition}

\begin{proof}
Take $\Omega_0\subset\subset\Omega$ and consider $\rho>0$ such that dist$(\Omega_0,\partial\Omega)>\rho$.
Let~$\eta\in C_0^\infty(B_1,[0,+\infty))$ be such that
$$ \int_{B_1}\eta(x)\,dx=1$$ and
$$\eta_\rho(x):=\frac1{\rho^{n}}\eta\left(\frac{x}{\rho}\right).$$
Define $v_s=u_s\ast\eta_\rho$ (extending~$u_s$ by zero outside~$\Omega$ for this convolution).

Thus, for every multiindex $\alpha$, 
$$D^\alpha v_s=u_s\ast D^\alpha \eta_\rho,$$ and then $$D^\alpha v_s\to u\ast D^\alpha \eta_\rho,$$ uniformly in $\overline{\Omega_0}$ as $s\nearrow1$. Hence 
\begin{equation}\label{convMoll}
\begin{split}
&v_s\to v:=u\ast \eta_\rho \;\mbox{ in }C^2(\overline{\Omega_0})\\
&{\mbox{and}}\qquad \|v_s\|_{C^2(\overline{\Omega_0})}\leq C(\rho,n){\|u_s\|_{L^2(\Omega_0+B_\rho)}}.\end{split}
\end{equation}
Moreover, by Jensen's inequality,
\begin{equation*}\begin{split}(v_s(x)-v_s(y))^2&=\bigg({\int_{\R^n}}\eta_\rho(z)(u_s(x-z)-u_s(y-z))\,dz\bigg)^2\\&\leq {\int_{\R^n}}\eta_\rho(z)(u_s(x-z)-u_s(y-z))^2\,dz,\end{split}\end{equation*}
and thus
\begin{equation}\begin{split}\label{ineqvs}&
\iint_{\Omega_0\times\Omega_0}\frac{(v_s(x)-v_s(y))^2}{|x-y|^{n+2s}} \,dx\,dy\\&\qquad\leq \iint_{\Omega_0\times\Omega_0}{\int_{\R^n}}\eta_\rho(z)\frac{(u_s(x-z)-u_s(y-z))^2}{|x-y|^{n+2s}}\,dz\,dx\,dy\\
&\qquad\leq {\int_{\R^n}}\eta_\rho(z)\iint_{\Omega\times\Omega}\frac{(u_s(x)-u_s(y))^2}{|x-y|^{n+2s}}\,dx\,dy\,dz \\&\qquad= \iint_{\Omega\times\Omega}\frac{(u_s(x)-u_s(y))^2}{|x-y|^{n+2s}}\,dx\,dy.
\end{split}\end{equation}

Let $\tilde{\Omega}_0\subset\subset\Omega_0$. Then, there exists {$\delta\in(0,1)$} such that $x+h\in \Omega_0$ for every $x\in\tilde{\Omega}_0$ and~$|h|<\delta$. We thus perform a Taylor expansion and arrive at
$$(v_s(x+h)-v_s(x))^2\geq (\nabla v_s\cdot h)^2-A^2|h|^3,$$
where {$A$ is a sufficiently large constant controlling the~$C^2$ norms of~$v_s$, uniformly for~$s$ close to~$1$, by~\eqref{convMoll} and local~$L^2$ convergence.} Therefore, we find that
\begin{equation}\begin{split}
\iint_{\Omega_0\times\Omega_0}\frac{(v_s(x)-v_s(y))^2}{|x-y|^{n+2s}} \,dx\,dy&\geq \int_{\tilde{\Omega}_0}\int_{\{|h|<\delta\}}\frac{(\nabla v_s\cdot h)^2-A^2|h|^3}{|h|^{n+2s}}\,dh\,dx\\
&{\ge}\frac{|\partial B_1|\delta^{2-2s}}{n(2-2s)}\int_{\tilde{\Omega}_0}|\nabla v_s|^2\,dx-|\partial B_1| A^2|\tilde{\Omega}_0|.
\end{split}\end{equation}
Therefore, substituting into~\eqref{ineqvs}, using that
$$\lim_{s\nearrow1}\frac{c_{n,s}|\partial B_1|}{2n(2-2s)}=1$$ and the fact that $\nabla v_s\to \nabla v$ uniformly {on} $\overline{\Omega_0}$, we get
$$\liminf_{s\nearrow1}\frac{c_{n,s}}{2}\iint_{\Omega\times\Omega}\frac{(u_s(x)-u_s(y))^2}{|x-y|^{n+2s}}\,dx\,dy\ \ge\int_{\tilde{\Omega}_0}|\nabla v|^2\,dx.$$

Sending now~$\rho\searrow0$ and exhausting~$\Omega_0$ by~$\tilde\Omega_0$, we get
$$\liminf_{s\nearrow1}\frac{c_{n,s}}{2}\iint_{\Omega\times\Omega}\frac{(u_s(x)-u_s(y))^2}{|x-y|^{n+2s}}\,dx\,dy\ \ge\int_{\Omega_0}|\nabla u|^2\,dx.$$
Finally, taking a sequence $\Omega_{0,k}\subset\subset\Omega$ such that $\Omega_{0,k}\to \Omega$ as $k\to+\infty$, we obtain the desired result.
\end{proof}

We now discuss the local compactness and convergence of the energy.

\begin{lemma}
%[Local compactness and convergence of energy forms]
\label{lem:energy-compactness}
Let~$s_j\in(0,1)$ with~$s_j\nearrow1$.
Suppose that either~$G_j=G=\R^n_+$, or~$G_j=\{x_n>-b_j\}$ with~$b_j\to+\infty$ and~$G=\R^n$.

For measurable functions~$h$, $k:\R^n\to\R$ with finite energy, set
\begin{equation}\label{compact:forms}
\begin{split}
{E_j}(h)&:=\frac{c_{n,s_j}}2
\iint_{\R^{2n}\setminus(G_j^c\times G_j^c)}
\frac{|h(x)-h(y)|^2}{|x-y|^{n+2s_j}}\,dx\,dy+\int_{G_j}h^2\,dx\\
{\mbox{and}}\qquad
{B_j}(h,k)&:=\frac14\big({E_j}(h+k)-{E_j}(h-k)\big).
\end{split}
\end{equation}
Also, for~$h,k\in H^1(G)$, set
\begin{equation}\label{compact:local-forms}
\begin{aligned}
E_G(h)&:=\int_G(|\nabla h|^2+h^2)\,dx,\\
B_G(h,k)&:=\int_G(\nabla h\cdot\nabla k+hk)\,dx.
\end{aligned}
\end{equation}
Let~$f_j:\R^n\to\R$ be measurable functions with~$E_j(f_j)\le C$ for every~$j$.

Then, up to a subsequence, there exists~$f\in H^1(G)$ such that
\begin{equation*}
f_j\to f\quad\text{in }L^2(B_R\cap G)
\qquad\text{for every }R>0
\end{equation*}
and
\begin{equation}\label{compact:lower}
E_G(f)\le\liminf_{j\to+\infty}E_j(f_j).
\end{equation}
Moreover, for every~$\Phi\in H^1(\R^n)$,
\begin{equation}\label{compact:test}
\lim_{j\to+\infty}E_j(\Phi)=E_G(\Phi)
\end{equation}
and
\begin{equation}\label{compact:mixed}
\lim_{j\to+\infty}B_j(f_j,\Phi)=B_G(f,\Phi).
\end{equation}
\end{lemma}

\begin{proof}
Let~$F_j$ be the even extension of~$f_j|_{G_j}$ across~$\partial G_j$. Let also~$y^*$ denote the reflection of~$y\in G_j$. In this way, for all~$x,y\in G_j$, we have that~$|x-y^*|\ge|x-y|$ and consequently
\begin{equation}\label{SK.92ed2rfuvhibrk}\int_{\R^n}(1+|\xi|^{2s_j})|\widehat F_j(\xi)|^2\,d\xi
\le4{E_j}(f_j)\le4C,\end{equation}
where~$\widehat F_j$ denotes the Fourier transform of~$F_j$.

Also, fix~$t\in(0,1)$. For large~$j$, we have that~$s_j\ge t$, so~\eqref{SK.92ed2rfuvhibrk} gives
\begin{equation}\label{compact:lower-order-bound}
\begin{split}&
\|F_j\|_{H^t(\R^n)}^2
=\int_{\R^n}(1+|\xi|^{2t})|\widehat F_j(\xi)|^2\,d\xi\\
&\qquad\le2\int_{\R^n}(1+|\xi|^{2s_j})|\widehat F_j(\xi)|^2\,d\xi
\le8E_j(f_j)\le8C.\end{split}
\end{equation}
Therefore, the weak compactness in~$L^2(\R^n)$ and the fractional Rellich--Kondrachov theorem~\cite[Theorem~7.1]{MR2944369}, applied in each ball, give a subsequence with~$F_j\rightharpoonup F$ in~$L^2(\R^n)$ and~$F_j\to F$ in~$L^2_{\rm loc}(\R^n)$. Thanks to a diagonal argument, we can suppose that the subsequence is independent of the ball.

Moreover, Plancherel's theorem gives that~$\widehat F_j\rightharpoonup\widehat F$ in~$L^2(\R^n)$. Also, in each ball~$B_R$, the weights~$1+|\xi|^{2s_j}$ converge uniformly to~$1+|\xi|^2$. The weak lower semicontinuity and~\eqref{SK.92ed2rfuvhibrk} therefore imply
\begin{equation}\label{compact:fourier-limit}
\int_{B_R}(1+|\xi|^2)|\widehat F(\xi)|^2\,d\xi
\le\liminf_{j\to+\infty}
\int_{B_R}(1+|\xi|^{2s_j})|\widehat F_j(\xi)|^2\,d\xi
\le4C.
\end{equation}
Letting~$R\to+\infty$ in~\eqref{compact:fourier-limit} gives that~$F\in H^1(\R^n)$.

Consequently,~$f:=F|_G$ belongs to~$H^1(G)$, and~$f_j=F_j$ in~$B_R\cap G$ for all sufficiently large~$j$.
Therefore, $f_j\to f$ in~$L^2(B_R\cap G)$ for every~$R>0$.

Now, in the case~$G_j=G=\R^n_+$, we apply Proposition~\ref{prop:liminf} on~$G$ and Fatou's lemma to obtain~\eqref{compact:lower}.
Also, in the case~$G=\R^n$, we use Proposition~\ref{prop:liminf} and Fatou's lemma in each~$B_R\subset G_j$ for large~$j$, and then let~$R\to+\infty$ to obtain~\eqref{compact:lower}.

Next, to prove~\eqref{compact:test}, let~$\Phi\in H^1(\R^n)$.
For~$G_j=G=\R^n_+$, equation~\eqref{halfspace:fixed-test}, with~$s=s_j$, gives~\eqref{compact:test} directly.

For~$b_j\to+\infty$ and~$G=\R^n$, fix~$R>0$. For large~$j$,
we have that~$B_R\subset G_j$, hence
\begin{equation*}
\begin{split}
&E_j(\Phi)\ge\frac{c_{n,s_j}}2\iint_{B_R\times B_R}
\frac{|\Phi(x)-\Phi(y)|^2}{|x-y|^{n+2s_j}}\,dx\,dy
+\int_{B_R}\Phi^2\,dx\\{\mbox{and}}\qquad &
\liminf_{j\to+\infty}E_j(\Phi)
\ge\int_{B_R}(|\nabla\Phi|^2+\Phi^2)\,dx,
\end{split}
\end{equation*}
by~\cite[Theorem~1]{MR4093822}. Letting~$R\to+\infty$ gives 
an inequality in~\eqref{compact:test}. Moreover,
\begin{equation*}
\limsup_{j\to+\infty}E_j(\Phi)
\le\lim_{j\to+\infty}\|\Phi\|_{H^{s_j}(\R^n)}^2
=\|\Phi\|_{H^1(\R^n)}^2=E_G(\Phi),
\end{equation*} which completes the proof of~\eqref{compact:test}.

{To prove~\eqref{compact:mixed}, we use the Cauchy--Schwarz inequality and~\eqref{compact:test} to find that
\begin{equation*}
|B_j(f_j,\Phi)|^2\le E_j(f_j)E_j(\Phi)\le C.
\end{equation*}
Consequently, along any subsequence, we can extract a further subsequence on which both bounded quantities~$E_j(f_j)$ and~$B_j(f_j,\Phi)$ converge. On this subsequence, let}
$$a:=\lim_{j\to+\infty} {E_j}(f_j)
\qquad{\mbox{and}}\qquad
b:=\lim_{j\to+\infty} {B_j}(f_j,\Phi).$$
Thus, for every~$\tau\in\mathbb R$, the quadratic expansion of~${E_j}$ and~\eqref{compact:lower} yield that
\begin{equation*}
\begin{aligned}
{E_G}(f+\tau\Phi)
\le \liminf_{j\to+\infty} {E_j}(f_j+\tau\Phi)
=a+2\tau b+\tau^2{E_G}(\Phi).
\end{aligned}
\end{equation*}
Therefore, expanding~$E_G(f+\tau\Phi)$ by~\eqref{compact:local-forms} and cancelling the terms in~$\tau^2$ gives
\begin{equation*}
a-{E_G}(f)
+2\tau\big[b-{B_G}(f,\Phi)\big]\ge0.
\end{equation*}
Hence, since this is valid for all~$\tau\in\R$, we get that~$
b={B_G}(f,\Phi)$.
Consequently, every convergent subsequence of~$B_j(f_j,\Phi)$ has limit~$B_G(f,\Phi)$, and boundedness proves~\eqref{compact:mixed}.
\end{proof}

We can now establish the desired compactness statement.

\begin{lemma}\label{nustonu1}
Let~$s_0$ be as in Lemma~\ref{3e4rfOSJL:LEM}, and let~$s_j\in(s_0,1)$ with~$s_j\nearrow1$. Let~$u_{s_j}\in \mathcal{M}_{s_j}$. 

Then, there exists $u_1\in \mathcal{M}_1$ such that, up to a subsequence 
and up to tangential translations, 
\begin{equation}\label{MCOSM:AS}u_{s_j}\to u_1\qquad {\mbox{ in $ L^2(\R^n_+)$ and in~$ L^{p+1}(\R^n_+)$ as $j\to+\infty$}}.\end{equation}

For positive minimizers centered as in~\eqref{REARR:in}, the full sequence converges to the positive centered normalized classical minimizer~$U$.

Moreover,
\begin{equation}\label{nu:convergence}
\lim_{s\nearrow1}\nu_s=\nu_1.
\end{equation}
\end{lemma}

\begin{proof}
We use the short notation~$u_j:=u_{s_j}$. By Theorem~\ref{ADDITIO:PR}, it suffices, after selecting a subsequence of constant sign, changing sign if necessary, and making tangential translations, to consider positive minimizers centered as in~\eqref{REARR:in}. Then, by Lemmata~\ref{SI-en-b} and~\ref{3e4rfOSJL:LEM},
\begin{equation*}
\|u_j\|_{s_j,+}^2=\nu_{s_j}\le C,\qquad
\|u_j\|_{L^{p+1}(\R^n_+)}=1\qquad{\mbox{and}}\qquad 0\le u_j\le M.
\end{equation*}
After passing to a subsequence, we may assume that~$\nu_{s_j}\to\nu$.

Moreover, if~$\widetilde U$ is the even extension of~$U$,
we have that~$\widetilde U$ is a competitor for~$\nu_s$, therefore~$
\nu_s\le\|\widetilde U\|_{s,+}^2$.
Therefore, applying~\eqref{compact:test} in Lemma~\ref{lem:energy-compactness} with~$G_j=G=\R^n_+$ and~$\Phi=\widetilde U$ gives
\begin{equation}\label{compact:upper}
\nu:=\lim_{j\to+\infty}\nu_{s_j}
\le\lim_{j\to+\infty}\|\widetilde U\|_{s_j,+}^2
=\int_{\R^n_+}(|\nabla U|^2+U^2)\,dx=\nu_1.\end{equation}

We next show that a given ball contains a uniformly positive amount of mass
of the ground state. More precisely, let~$F_j$ be the even reflection of~$u_j$ across~$\partial\R^n_+$.
Let~$t\in(0,1)$ with~$p+1<2_t^*$ and~$r:=2+4t/n$.
Let also~$Q_k:=k+(0,1)^n$, with~$k\in\mathbb Z^n$.
We consider the
translates~$\eta_k$ of a fixed smooth cutoff. Namely, we
take~$\eta_k\in C^\infty_0(k+(-1,2)^n,[0,1])$ such that~$
\eta_k=1$ in~$Q_k$. 

We observe that
$$\sum_{k\in\Z^n}|\eta_k(x)-\eta_k(y)|^2\le C\min\{1,|x-y|^2\}$$
and therefore
$$\sum_{k\in\Z^n}\|\eta_kF_j\|_{H^t(\R^n)}^2\le C\|F_j\|_{H^t(\R^n)}^2.$$
Accordingly, by H\"older's inequality
and the Sobolev inequality, we find that
\begin{equation}\label{EFFEj}\begin{split}&
\|F_j\|_{L^r(\R^n)}^r\le
\left(\sup_{k\in\Z^n}\int_{Q_k}F_j^2\,dx\right)^{2t/n}
\sum_{k\in\Z^n}\|\eta_kF_j\|_{2_t^*}^2
\\&\qquad\le C\left(\sup_{k\in\Z^n}\int_{Q_k}F_j^2\,dx\right)^{2t/n}\|F_j\|_{H^t(\R^n)}^2.\end{split}\end{equation}

Next, we claim that
\begin{equation}\label{LIMASSa}
\liminf_{j\to+\infty}\sup_{k\in\Z^n}\int_{Q_k}F_j^2\,dx>0.\end{equation}
Indeed, suppose for the sake of contradiction, that, up to a subsequence,
$$\lim_{j\to+\infty}\sup_{k\in\Z^n}\int_{Q_k}F_j^2\,dx=0.$$
Then, by~\eqref{EFFEj},
$$\lim_{j\to+\infty}\|F_j\|_{L^r(\R^n)}=0.$$
{If~$2<p+1\le r$, by interpolation and recalling~\eqref{SK.92ed2rfuvhibrk} we get
\begin{equation*}
\lim_{j\to+\infty}\|F_j\|_{L^{p+1}(\R^n)}
\le \lim_{j\to+\infty}\|F_j\|_{L^2(\R^n)}^{\frac{2(r-p-1)}{(p+1)(r-2)}}
\|F_j\|_{L^r(\R^n)}^{\frac{r(p-1)}{(p+1)(r-2)}}
=0.
\end{equation*}
If instead~$p+1>r$, Lemma~\ref{3e4rfOSJL:LEM} gives that~$\|F_j\|_{L^\infty(\R^n)}\le M$, and consequently
\begin{equation*} \lim_{j\to+\infty}
\|F_j\|_{L^{p+1}(\R^n)}^{p+1}
\le M^{p+1-r} \lim_{j\to+\infty}\|F_j\|_{L^r(\R^n)}^r
=0.
\end{equation*}}
In any case, we obtain a contradiction with the fact that~$\|F_j\|_{L^{p+1}(\R^n)}^{p+1}=2$
and the proof of~\eqref{LIMASSa} is thereby complete.

As a byproduct of~\eqref{LIMASSa}, we can find~$R$, $\delta>0$ and~$y_j=(y'_j,y_{j,n})\in\overline{\R^n_+}$ such that
\begin{equation}\label{compact:mass}
\int_{B_R(y_j)\cap\R^n_+}u_j^2\,dx\ge\delta.
\end{equation}

We claim that
\begin{equation}\label{compact:mass:00b}
{\mbox{the sequence~$y'_j$ is bounded. }}\end{equation}

Indeed, suppose, by contradiction, that~$|y'_j|\to+\infty$ along a subsequence, and pick an integer~$N\ge2$.
Also, for~$k=0,\ldots,N-1$, define
\begin{equation}\label{compact:mass:centers}
 y_{j,k}:=|y'_j|\left(\cos\left(\frac{2\pi k}{N}\right)\,e_1
 +\sin\left(\frac{2\pi k}{N}\right)\,e_2\right)+y_{j,n}e_n.
\end{equation}
In this setting, the tangential rotational symmetry in~\eqref{REARR:in} and~\eqref{compact:mass} give that
\begin{equation}\label{compact:mass:rotations}
\int_{B_R(y_{j,k})\cap\R^n_+}u_j^2\,dx
=\int_{B_R(y_j)\cap\R^n_+}u_j^2\,dx\ge\delta.
\end{equation}
Besides, for distinct~$k,\ell\in\{0,\ldots,N-1\}$, the centers in~\eqref{compact:mass:centers} satisfy, for all sufficiently large~$j$,
\begin{equation*}
|y_{j,k}-y_{j,\ell}|
=2|y'_j|\left|\sin\frac{\pi(k-\ell)}N\right|
\ge2|y'_j|\sin\frac\pi N>2R.
\end{equation*}
In view of this observation, the balls~$B_R(y_{j,k})$ are pairwise disjoint, and~\eqref{compact:mass:rotations}, together with Lemma~\ref{SI-en-b}, yields that
\begin{equation*}
\begin{aligned}
N\delta\le\sum_{k=0}^{N-1}
\int_{B_R(y_{j,k})\cap\R^n_+}u_j^2\,dx\le\|u_j\|_{L^2(\R^n_+)}^2\le\nu_{s_j}\le C.
\end{aligned}
\end{equation*}
This produces a contradiction by choosing~$N>C/\delta$. The proof of~\eqref{compact:mass:00b} is thereby complete.

In fact, we have that
\begin{equation}\label{compact:mass:00c}
{\mbox{the sequence~$y_j$ is bounded. }}\end{equation}
To check this,
suppose not. Then, along a subsequence, 
we infer from~\eqref{compact:mass:00b}
that~$y_{j,n}\to+\infty$. We let~$v_j(x):=u_j(x+y_j)$ and set~$G_j:=\{x_n>-y_{j,n}\}$. Changing variables~$X:=x+y_j$ and~$Y:=y+y_j$ in~\eqref{compact:forms} gives
\begin{equation*}
E_j(v_j)=\|u_j\|_{s_j,+}^2=\nu_{s_j}\le C\qquad{\mbox{and}}
\qquad \int_{G_j}v_j^{p+1}\,dx=1.
\end{equation*}
Therefore, Lemma~\ref{lem:energy-compactness}, applied with~$f_j=:v_j$ and~$b_j:=y_{j,n}\to+\infty$, gives, up to a subsequence, a function~$v\in H^1(\R^n)$ such that~$v_j\to v$ in~$L^2(B_L)$ for every~$L>0$. 

Also, for large~$j$, a given ball~$B_R$ lies in~$G_j$, hence~\eqref{compact:mass} yields that
$$\int_{B_R}v_j^2\,dx\ge\delta.$$ Thus, passing to the limit,
$$\int_{B_R}v^2\,dx\ge\delta,$$ and therefore~$v$ does not vanish identically. 

Besides, for each~$L>0$, Fatou's lemma gives
that$$\int_{B_L}v^{p+1}\,dx\le1$$
and thus, letting~$L\to+\infty$,
$$\|v\|_{L^{p+1}(\R^n)}\le1.$$

Moreover, the translated weak equation in this setting becomes, for all~$\Phi\in C^\infty_0(\R^n)$, 
\begin{equation*}
B_j(v_j,\Phi)=\nu_{s_j}\int_{G_j}v_j^p\Phi\,dx.
\end{equation*}
Hence, by~\eqref{compact:mixed}, 
\begin{equation}\label{POThTH}\int_{\R^n}(\nabla v\cdot\nabla\Phi+v\Phi)\,dx=
\lim_{j\to+\infty}B_j(v_j,\Phi)=\nu\lim_{j\to+\infty}\int_{G_j}v_j^p\Phi\,dx.
\end{equation}

Additionally,
on the compact support of~$\Phi$
(which lies in~$G_j$ for large~$j$) the bounds~$
0\le v_j,v\le M$ give that
\begin{equation}\label{Pvjpoiuj}
|v_j^p-v^p|\le pM^{p-1}|v_j-v|.
\end{equation}
Thus, by the convergence in~$L^2_{\rm loc}(\R^n)$,
$$\nu\lim_{j\to+\infty}\int_{G_j}v_j^p\Phi\,dx=
{\nu}\int_{\R^n}v^p\Phi\,dx,$$
which, in tandem with~\eqref{POThTH}, establishes that, for all~$\Phi\in C^\infty_0(\R^n)$ (and for all~$\Phi\in H^1(\R^n)$ by density),
\begin{equation*}\int_{\R^n}(\nabla v\cdot\nabla\Phi+v\Phi)\,dx
=\nu\int_{\R^n}v^p\Phi\,dx.\end{equation*}

Taking~$\Phi:=v$ and setting~$a:=\|v\|_{L^{p+1}(\R^n)}\in(0,1]$, {the definition of~$\nu_1^\star$ in~\eqref{nu1:whole-space} gives}
\begin{equation*}
{\nu_1^\star}a^2\le\|v\|_{H^1(\R^n)}^2=\nu a^{p+1}\qquad{\mbox{and}}
\qquad {\nu_1^\star}\le\nu a^{p-1}\le\nu\le\nu_1.
\end{equation*}
This contradicts the fact that~$\nu_1<{\nu_1^\star}$, see Remark~\ref{relationnu}. Thus, the claim in~\eqref{compact:mass:00c} is established.

Next, we apply Lemma~\ref{lem:energy-compactness} with~$G_j=G=\R^n_+$ and~$f_j=u_j$. In this way,
up to a subsequence, there exists~$u\in H^1(\R^n_+)$ such that~$u_j\to u$ in~$L^2(B_L\cap\R^n_+)$ for every~$L>0$. 

On account of~\eqref{compact:mass:00c}, we can also find~$L_0>0$ such that~$B_R(y_j)\subset B_{L_0}$ for all~$j\in\N$. 
Therefore, we deduce from~\eqref{compact:mass} that$$\int_{B_{L_0}\cap\R^n_+}u^2\,dx\ge\delta,$$ showing that~$u$ does not vanish identically. 

We now consider the weak equation associated with~$u_j$
and pass to the limit as~$j\to+\infty$. Namely, given~$\Phi\in C^\infty_0(\R^n)$,
starting with
\begin{equation*}
B_j(u_j,\Phi)=\nu_{s_j}\int_{\R^n_+}u_j^p\Phi\,dx,
\end{equation*}
using~\eqref{compact:mixed} and~\eqref{Pvjpoiuj}
(here for~$u_j$ and~$u$ instead of~$v_j$ and~$v$), and passing to the limit, we obtain that
\begin{equation}\label{compact:limit-equation}
\int_{\R^n_+}(\nabla u\cdot\nabla\Phi+u\Phi)\,dx
=\nu\int_{\R^n_+}u^p\Phi\,dx.
\end{equation}
By density, we are allowed to take~$\Phi:=u$ in~\eqref{compact:limit-equation}. Then,
Fatou's lemma gives that$$a:=\|u\|_{L^{p+1}(\R^n_+)}\in(0,1].$$ Consequently,
\begin{equation*}
\nu_1a^2\le\|u\|_{H^1(\R^n_+)}^2=\nu a^{p+1}
\qquad{\mbox{and}}\qquad \nu_1\le\nu a^{p-1}\le\nu\le\nu_1.
\end{equation*}
Therefore
\begin{equation}\label{Ajmda1}
{\mbox{$a=1$, $\nu=\nu_1$ and~$u\in\mathcal M_1$.}}
\end{equation}

Additionally, Lemma~\ref{3e4rfOSJL:LEM} gives that~$0\le u_j\le M$ in~$\R^n_+$, and the same bound holds for~$u$. Consequently, for every~$R>0$,
\begin{equation*}\lim_{j\to+\infty}
\int_{B_R\cap\R^n_+}|u_j-u|^{p+1}\,dx
\le(2M)^{p-1}\lim_{j\to+\infty}\int_{B_R\cap\R^n_+}|u_j-u|^2\,dx=0.
\end{equation*}

Hence, using the normalization in~\eqref{MIN:GRO:ST} 
and the fact that~$\|u\|_{L^{p+1}(\R^n_+)}=1$ (due to~\eqref{Ajmda1}),
we obtain, for every~$R>0$,
\begin{equation}\label{KP:SD:k2p1WT}
\begin{split}&
\lim_{j\to+\infty}\int_{\R^n_+\setminus B_R}u_j^{p+1}\,dx
=1-\lim_{j\to+\infty}\int_{B_R\cap\R^n_+}u_j^{p+1}\,dx\\
&\qquad=1-\int_{B_R\cap\R^n_+}u^{p+1}\,dx
=\int_{\R^n_+\setminus B_R}u^{p+1}\,dx.
\end{split}
\end{equation}
Thus, since
$$|u_j-u|^{p+1}\le2^{p}(u_j^{p+1}+u^{p+1}),$$
we deduce from~\eqref{KP:SD:k2p1WT},
by first letting~$j\to+\infty$ and then~$R\to+\infty$, 
that~$u_j\to u$ in~$L^{p+1}(\R^n_+)$.

To complete the proof of~\eqref{MCOSM:AS}, it remains to
prove the convergence in~$L^2(\R^n_+)$.
With this objective in mind, we use Proposition~\ref{prop:liminf} to get that$$\liminf_{j\to+\infty}\nu_{s_j}-\|u_j\|_{L^2(\R^n_+)}^2\ge\|\nabla u\|_2^2.$$ Hence,
\begin{equation*}
\limsup_{j\to+\infty}\|u_j\|_2^2\le\nu_1-\|\nabla u\|_2^2=\|u\|_2^2.
\end{equation*}
Together with the weak convergence in~$L^2(\R^n_+)$, this leads to strong convergence in~$L^2(\R^n_+)$, thus finishing the proof of~\eqref{MCOSM:AS}.

Finally, for positive centered minimizers, thanks to Remark~\ref{relationnu} and classical uniqueness we can conclude that~$u$ is the positive centered classical minimizer~$U$ (see~\cite[Theorem~1.1]{FKNOTES}). Consequently, every subsequence of~$u_{s_j}$ has a further subsequence with~$u_{s_j}\to U$ as in~\eqref{MCOSM:AS} and~$\nu_{s_j}\to\nu_1$ by~\eqref{Ajmda1}. The uniqueness of both subsequential limits proves the convergence of the full sequence and~\eqref{nu:convergence}.
\end{proof}

 \section{Nondegeneracy and proof of Theorem~\ref{th:nondegeneracy}}\label{nondegeneracy:section}
 The purpose of this section is to prove Theorem~\ref{th:nondegeneracy}, which will be obtained as a consequence of two facts: the asymptotics proven in
Section~\ref{asymptotics} and a precise description of the form of $w$. 

We use the energy space~$H^s_+$ throughout this section. A function~$w\in H^s_+$ is a weak solution of~\eqref{lin_problem} if, for all~$\varphi\in H^s_+$,
\begin{equation}\label{WeakForm_eq}
(w,\varphi)_{s,+}
=p\nu_s\int_{\R^n_+}u_s^{p-1}w\varphi\,dx.
\end{equation}
The bilinear form~$(\cdot,\cdot)_{s,+}$ is the one defined in footnote~\ref{ihqkdncs:MAsxRT}.

We now present some auxiliary lemmata: 
\begin{lemma}\label{lemmaJs}
Let $\varphi\in H^s_+$ such that $(u_s,\varphi)_{s,+}=0$. Then
\begin{equation}\label{lemmaJs:form}
0\leq \mathcal{J}_s''(u_s,\nu_s)[\varphi,\varphi]=\|\varphi\|_{s,+}^2-p\nu_s\int_{\R^n_+}|u_s|^{p-1}\varphi^2\,dx.
\end{equation}
 \end{lemma}

The proof of Lemma~\ref{lemmaJs} follows in the same way as \cite[Lemma~3.1]{Fall-Valdinoci} just by noticing that~$(u_s,\varphi)_{s,+}=0$ implies
 $$\int_{\R^n_+}{u_s^p}\varphi\,dx=0,$$
simply testing with $\varphi$ in the weak formulation of~\eqref{problem_nu}.

We also have:

\begin{lemma}\label{lemmaA1}
Let~$u_s$ be centered as in~\eqref{REARR:in}.
%, using the tangential translation supplied by Theorem~\ref{ADDITIO:PR}.
Let~$i\in\{1,\ldots,n-1\}$ and let~$g^i\in H^s_+$ be odd in~$x_i$. Define
\begin{equation*}
\psi(x):=\begin{cases}
|g^i(x)|,& {\mbox{ if }}x_i>0,\\
-|g^i(x)|,&{\mbox{ if }}x_i<0,\\
0,& {\mbox{ if }}x_i=0.
\end{cases}
\end{equation*}
Then~$\psi\in H^s_+$ and
\begin{equation}\label{A1}
0\le\mathcal J_s''(u_s,\nu_s)[\psi,\psi]
\le\mathcal J_s''(u_s,\nu_s)[g^i,g^i].
\end{equation}

Moreover, if~$g^i$ also solves~\eqref{lin_problem}, then it has one sign almost everywhere in~$H_i:=\{x_i>0\}$, unless it vanishes identically.
\end{lemma}

\begin{proof} For all~$x,y\in H_i$ with~$x\ne y$, we have that~$|x-(y-2y_i e_i)|^2-|x-y|^2=4x_i y_i>0$, therefore
\begin{equation}\label{A1:kernel-positive}
\frac1{|x-y|^{n+2s}}
-\frac1{|x-(y-2y_i e_i)|^{n+2s}}>0.
\end{equation}
As a result, 
\begin{eqnarray*}
&&\frac{(g^i(x)-g^i(y))^2-(|g^i(x)|-|g^i(y)|)^2}
{|x-y|^{n+2s}}+\frac{(g^i(x)+g^i(y))^2-(|g^i(x)|+|g^i(y)|)^2}
{|x-(y-2y_i e_i)|^{n+2s}}\\
&&\qquad=2 \Big(|g^i(x)g^i(y)|-g^i(x)g^i(y)\Big)
\left(
\frac1{|x-y|^{n+2s}}-\frac1{|x-(y-2y_i e_i)|^{n+2s}}
\right)\ge0.\end{eqnarray*}
Therefore,
splitting the integrals according to whether each variable belongs to~$H_i$ or its complement and using odd symmetry, we see that
\begin{equation}\label{A1:reflection}
\begin{split}
\|\psi\|_{s,+}^2
&=c_{n,s}\iint_{(H_i\times H_i)\setminus(\R^n_-\times\R^n_-)}
\left[
\frac{(|g^i(x)|-|g^i(y)|)^2}{|x-y|^{n+2s}}+\frac{(|g^i(x)|+|g^i(y)|)^2}{|x-(y-2y_i e_i)|^{n+2s}}
\right]\,dx\,dy\\&\qquad\qquad
+2\int_{H_i\cap\R^n_+}(g^i)^2\,dx\\
&\le c_{n,s}\iint_{(H_i\times H_i)\setminus(\R^n_-\times\R^n_-)}
\left[
\frac{(g^i(x)-g^i(y))^2}{|x-y|^{n+2s}}
+\frac{(g^i(x)+g^i(y))^2}{|x-(y-2y_i e_i)|^{n+2s}}
\right]\,dx\,dy\\&\qquad\qquad
+2\int_{H_i\cap\R^n_+}(g^i)^2\,dx\\
&=\frac{c_{n,s}}2
\iint_{\R^{2n}\setminus(\R^n_-\times\R^n_-)}
\frac{(g^i(x)-g^i(y))^2}{|x-y|^{n+2s}}\,dx\,dy
+\int_{\R^n_+}(g^i)^2\,dx\\
&=\|g^i\|_{s,+}^2,
\end{split}
\end{equation}
which is finite, proving that~$\psi\in H^s_+$.

Also, if either~$\varphi:=g^i$ or~$\varphi:=\psi$,
we have that~$\varphi\in H^s_+$ and it is odd in~$x_i$.
Therefore, testing~\eqref{problem_nu} with~$\varphi$, changing variables~$x\mapsto x-2x_i e_i$ in~$\R^n_+\setminus H_i$, and using~\eqref{REARR:in}, we get
\begin{equation*}
\begin{split}
(u_s,\varphi)_{s,+}
&=\nu_s\int_{\R^n_+}u_s^p(x)\varphi(x)\,dx\\
&=\nu_s\int_{H_i\cap\R^n_+}
\Big(u_s^p(x)\varphi(x)
+u_s^p(x-2x_i e_i)\varphi(x-2x_i e_i)\Big)\,dx\\
&=\nu_s\int_{H_i\cap\R^n_+}
u_s^p(x)\big(\varphi(x)-\varphi(x)\big)\,dx\\&=0.
\end{split}
\end{equation*}
In light of this, we can exploit Lemma~\ref{lemmaJs} with~$\varphi$.
In this way, recalling also~\eqref{A1:reflection}, we obtain that
\begin{equation}\label{A1:form-comparison}
\begin{aligned}
0\le\mathcal J_s''(u_s,\nu_s)[\psi,\psi]
&=\|\psi\|_{s,+}^2
-p\nu_s\int_{\R^n_+}u_s^{p-1}\psi^2\,dx\\
&=\|\psi\|_{s,+}^2
-p\nu_s\int_{\R^n_+}u_s^{p-1}(g^i)^2\,dx\\
&\le\|g^i\|_{s,+}^2
-p\nu_s\int_{\R^n_+}u_s^{p-1}(g^i)^2\,dx\\
&=\mathcal J_s''(u_s,\nu_s)[g^i,g^i],
\end{aligned}
\end{equation}
which proves~\eqref{A1}.

Next, suppose that~$g^i$ is a weak solution of~\eqref{lin_problem}.
Therefore, taking~$w=\varphi=g^i$ in~\eqref{WeakForm_eq}, we see that
\begin{equation}\label{eqGi3}
\mathcal J_s''(u_s,\nu_s)[g^i,g^i]
=\|g^i\|_{s,+}^2-p\nu_s\int_{\R^n_+}u_s^{p-1}(g^i)^2\,dx=0.
\end{equation}

We observe that if~$g^i=0$ almost everywhere in~$\R^n_+$, then~\eqref{eqGi3} yields
\begin{equation*}
\|g^i\|_{s,+}^2
=p\nu_s\int_{\R^n_+}u_s^{p-1}(g^i)^2\,dx=0.
\end{equation*} Consequently,
\begin{equation*}
0=\|g^i\|_{s,+}^2
=c_{n,s}\int_{\R^n_-}(g^i(y))^2
\left(\int_{\R^n_+}\frac{dx}{|x-y|^{n+2s}}\right)\,dy,
\end{equation*} which entails that~$g^i=0$ almost everywhere in~$\R^n_-$, and hence in~$\R^n$.

If instead~$g^i$ does not vanish identically in~$\R^n_+$, then by odd symmetry in~$x_i$,
\begin{equation*}
0<\int_{H_i\cap\R^n_+}(g^i)^2\,dx
=\frac12\int_{\R^n_+}(g^i)^2\,dx<+\infty.
\end{equation*}
Accordingly, we can find
\begin{equation}\label{A1:nonzero-halfspace}
{\mbox{$x\in H_i\cap\R^n_+$ with~$g^i(x)\ne0$.}}
\end{equation}

Moreover, from~\eqref{A1:form-comparison} and~\eqref{eqGi3} we infer
that
$$0\le\mathcal J_s''(u_s,\nu_s)[\psi,\psi]
\le\mathcal J_s''(u_s,\nu_s)[g^i,g^i]=0$$
and \begin{equation*}
\begin{split}&
\|g^i\|_{s,+}^2-\|\psi\|_{s,+}^2
=\mathcal J_s''(u_s,\nu_s)[g^i,g^i]
-\mathcal J_s''(u_s,\nu_s)[\psi,\psi]+p\nu_s\int_{\R^n_+}u_s^{p-1}
\big((g^i)^2-\psi^2\big)\,dx=0.
\end{split}
\end{equation*}
Consequently, recalling~\eqref{A1:reflection}, 
\begin{eqnarray*}
&&\iint_{(H_i\times H_i)\setminus(\R^n_-\times\R^n_-)}
\big(|g^i(x)g^i(y)|-g^i(x)g^i(y)\big)\left(
\frac1{|x-y|^{n+2s}}
-\frac1{|x-(y-2y_i e_i)|^{n+2s}}
\right)\,dx\,dy\\
&&\quad=\frac{\|g^i\|_{s,+}^2-\|\psi\|_{s,+}^2}{2c_{n,s}}=0.
\end{eqnarray*}
From this and~\eqref{A1:kernel-positive} we obtain that, for almost every~$(x,y)\in(H_i\times H_i)\setminus(\R^n_-\times\R^n_-)$,
\begin{equation*}
|g^i(x)g^i(y)|-g^i(x)g^i(y)=0,
\end{equation*} hence~$g^i(x)g^i(y)\ge0$.

This observation, together with~\eqref{A1:nonzero-halfspace},
entails that~$g^i(x)g^i(y)\ge0$ for almost every~$y\in H_i$.
Therefore, $g^i$ does not change sign almost everywhere in~$H_i$.
\end{proof}

Now we can prove the result that will be the core of the proof of Theorem~\ref{th:nondegeneracy}. Let us emphasize that, even if we will only prove the nondegeneracy for~$s$ close to 1, the following decomposition holds true for every~$s\in (0,1)$. 
 
 \begin{proposition}\label{decompw}
 Let $s\in (0,1)$ and let~$w\in H^s_+$ be a weak solution of~\eqref{lin_problem}, with~$u_s$ centered as in~\eqref{REARR:in}. Then
\begin{equation}\label{W:LS}
w(x',x_n)=f_0(|x'|,x_n)+\sum_{i=1}^{n-1}c_i\partial_{x_i}u_s,\end{equation}
 where~$c_i\in\R$ and, for all~$r>0$, we have that
\begin{equation}\label{W:LS2}f_0(r, x_n):=\frac1{|\mathbb{S}^{n-2}|}\int_{\mathbb{S}^{n-2}}w(r\theta,x_n)\,d\theta .\end{equation}
 \end{proposition}

\begin{proof}
Let~$i\in\{1,\ldots,n-1\}$. We claim that
\begin{equation}\label{0qoejld23:TSBqP}
{\mbox{$\partial_{x_i}u_s$ belongs to~$H^s_+$ and solves~\eqref{lin_problem} weakly.}}\end{equation} For this, given~$h\in\R\setminus\{0\}$, we subtract the weak equation~\eqref{problem_nu} from its tangential translate and test with~$u_s(\cdot+he_i)-u_s$. Hence, using~\eqref{OFCLA}, we obtain that
\begin{equation}\label{oqiqe30TBC} \|u_s(\cdot+he_i)-u_s\|_{s,+}^2
\le p\nu_s\|u_s\|_{L^\infty(\R^n_+)}^{p-1}
\|u_s(\cdot+he_i)-u_s\|_{L^2(\R^n_+)}^2.\end{equation}

Now, let~$\widetilde u_s(x',x_n):=u_s(x',|x_n|)$ be the even extension of~$u_s|_{\R^n_+}$.
Since~$i<n$, tangential a translation commutes with the even reflection.
Therefore, applying~\eqref{SK.92ed2rfuvhibrk} to~$u_s(\cdot+he_i)-u_s$, we find that\footnote{For better readability, here we are using~${\mathcal{F}}u$ to denote the Fourier transform of~$u$.}
\begin{equation}\label{decomp:fourier-difference}
\begin{split}
&\int_{\R^n}|\xi|^{2s}|e^{ih\xi_i}-1|^2
|\mathcal{F}\,{\widetilde u_s}(\xi)|^2\,d\xi=\int_{\R^n}|\xi|^{2s}
\left|{\mathcal{F}}(\widetilde u_s(\cdot+he_i)-\widetilde u_s)(\xi)\right|^2\,d\xi\\
&\qquad\le4\|u_s(\cdot+he_i)-u_s\|_{s,+}^2.
\end{split}
\end{equation}

Now let~$R\ge1$, independent of~$h$, such that
\begin{equation}\label{cjnkewgtiewfgbchksa45678ghj}
R^{2s}\ge4p\nu_s\|u_s\|_{L^\infty(\R^n_+)}^{p-1}.\end{equation}
Thus, Plancherel's theorem and~\eqref{decomp:fourier-difference} yield
\begin{equation*}
\begin{split}
&\|u_s(\cdot+he_i)-u_s\|_{L^2(\R^n_+)}^2=\frac12\int_{\R^n}|e^{ih\xi_i}-1|^2
|\mathcal{F}\,{\widetilde u_s}(\xi)|^2\,d\xi\\
&\quad\le\frac{R^2h^2}2
\int_{\{|\xi|\le R\}}|\mathcal{F}\,{\widetilde u_s}(\xi)|^2\,d\xi+\frac1{2R^{2s}}
\int_{\{|\xi|>R\}}|\xi|^{2s}|e^{ih\xi_i}-1|^2
|\mathcal{F}\,{\widetilde u_s}(\xi)|^2\,d\xi\\
&\quad\le R^2h^2\|u_s\|_{L^2(\R^n_+)}^2
+\frac2{R^{2s}}\|u_s(\cdot+he_i)-u_s\|_{s,+}^2.
\end{split}
\end{equation*}
Consequently, recalling~\eqref{cjnkewgtiewfgbchksa45678ghj}, we obtain that
\begin{equation*}\begin{split}&
p\nu_s  \|u_s\|_{L^\infty(\R^n_+)}^{p-1}
\|u_s(\cdot+he_i)-u_s\|_{L^2(\R^n_+)}^2\\
&\quad\le p\nu_s R^2h^2 \|u_s\|_{L^\infty(\R^n_+)}^{p-1}\|u_s\|_{L^2(\R^n_+)}^2
+\frac{2p\nu_s}{R^{2s}}\|u_s\|_{L^\infty(\R^n_+)}^{p-1}\|u_s(\cdot+he_i)-u_s\|_{s,+}^2\\
&\quad\le Ch^2\|u_s\|_{L^2(\R^n_+)}^2
+\frac{1}{2}\|u_s(\cdot+he_i)-u_s\|_{s,+}^2,
\end{split}\end{equation*}
with~$C:=p\nu_s\|u_s\|_{L^\infty(\R^n_+)}^{p-1}R^2$ independent of~$h$.

Combining this with~\eqref{oqiqe30TBC}, we find that
$$  \frac1{h^2}\|u_s(\cdot+he_i)-u_s\|_{s,+}^2 \le 2C \|u_s\|_{L^2(\R^n_+)}^2
.$$
Sending~$h\to0$, we thereby conclude that~$\partial_{x_i}u_s\in H^s_+$. 

Besides, passing to the limit in the difference of the two weak equations of~\eqref{problem_nu}, divided by~$h$, we also show that~$\partial_{x_i}u_s$ solves~\eqref{lin_problem} weakly. This completes the proof of~\eqref{0qoejld23:TSBqP}.

By~\eqref{REARR:in} and~\eqref{eqRedde}, we also have
\begin{equation}\label{UYUYi8}
{\mbox{$\partial_{x_i}u_s$ is odd in~$x_i$ and strictly negative on~$H_i\cap\R^n_+$.}}
\end{equation}

We now identify all solutions of~\eqref{lin_problem}
which are odd in~$x_i$. We claim that
\begin{equation}\label{decomp:odd}
\text{if $g^i\in H^s_+$ is a weak solution of~\eqref{lin_problem}
that is odd in $x_i$, then $
g^i=c_i\partial_{x_i}u_s$,}
\end{equation}
for a suitable~$c_i\in\R$.
To prove this, let~$g^i$ be such a solution. By~\eqref{0qoejld23:TSBqP} and~\eqref{UYUYi8}, we have that
\begin{equation*}\int_{H_i\cap\R^n_+}(\partial_{x_i}u_s)^2\,dx\in(0,+\infty)\end{equation*}
and, by the Cauchy--Schwarz inequality, 
$$\int_{H_i\cap\R^n_+}|g^i|\,|\partial_{x_i}u_s|\,dx\in{[0,+\infty)}.$$
These observations allow us to define
\begin{equation*}
c_i:=\frac{\displaystyle\int_{H_i\cap\R^n_+}g^i\partial_{x_i}u_s\,dx}
{\displaystyle\int_{H_i\cap\R^n_+}(\partial_{x_i}u_s)^2\,dx}.
\end{equation*}
Equivalently,
\begin{equation}\label{decomp:orthogonality}
\int_{H_i\cap\R^n_+}
(g^i-c_i\partial_{x_i}u_s)\partial_{x_i}u_s\,dx=0.
\end{equation}

Now, by~\eqref{0qoejld23:TSBqP}, \eqref{UYUYi8} and linearity of~\eqref{lin_problem}, the difference~$g^i-c_i\partial_{x_i}u_s$ is an odd weak solution of~\eqref{lin_problem} in~$H^s_+$. Hence, 
if this difference does not vanish identically in~$H_i$,
in virtue of Lemma~\ref{lemmaA1},
it has one sign in~$H_i$. 
This fact, together with the strict negativity in~\eqref{UYUYi8} and~\eqref{decomp:orthogonality}, entails that~$g^i-c_i\partial_{x_i}u_s=0$ almost everywhere in~$H_i\cap\R^n_+$, and therefore in~$\R^n_+$ by oddsymmetry, and then in~$\R^n$ thanks to
the Neumann condition in~\eqref{lin_problem}. This proves~\eqref{decomp:odd}.

In addition, in light of~\eqref{REARR:in}, we know that~$\partial_{x_i}u_s$ is even in every tangential coordinate different from~$x_i$.
Therefore, by~\eqref{0qoejld23:TSBqP} and~\eqref{UYUYi8} we have that, for all~$i,j\in\{1,\ldots,n-1\}$,
\begin{equation}\label{decomp:gradient-orthogonality}
\begin{split}&
\int_{\R^n_+}\partial_{x_i}u_s\,\partial_{x_j}u_s\,dx=0
\quad\text{if }i\ne j,\\ {\mbox{and }}\qquad&
0<\int_{\R^n_+}(\partial_{x_i}u_s)^2\,dx
=\int_{\R^n_+}(\partial_{x_1}u_s)^2\,dx<+\infty.
\end{split}
\end{equation}
On this account, for the given weak solution~$w$ of~\eqref{lin_problem}, the coefficients
\begin{equation}\label{decomp:coefficients}
c_i:=\frac{\displaystyle\int_{\R^n_+}w\,\partial_{x_i}u_s\,dx}
{\displaystyle\int_{\R^n_+}(\partial_{x_i}u_s)^2\,dx},
\qquad i=1,\ldots,n-1,
\end{equation}
are well-defined by the Cauchy--Schwarz inequality and~\eqref{decomp:gradient-orthogonality}.

Now, let~$\theta\in\mathbb S^{n-2}$ and~$\varphi\in H^s_+$
and notice that the reflection~$x\mapsto x-2(x'\cdot\theta)(\theta,0)$ preserves distances and~$\R^n_\pm$.
Also, by~\eqref{REARR:in}, $u_s$ is invariant under this reflection.
Therefore, we set
$$
\tilde\varphi(x):=\varphi(x-2(x'\cdot\theta)(\theta,0))\qquad
{\mbox{and}}\qquad\tilde w(x):=w(x-2(x'\cdot\theta)(\theta,0))$$
and,
changing variables under the above reflection in the norm and in~\eqref{WeakForm_eq}, we find that
\begin{equation}\label{decomp:reflection-invariance}
\begin{aligned}
&\big\|\tilde\varphi\big\|_{s,+}
=\|\varphi\|_{s,+}\\{\mbox{and }}\quad
&\mathcal J_s''(u_s,\nu_s)
\big[\tilde w,\,\varphi\big]=\mathcal J_s''(u_s,\nu_s)
\big[w,\,\tilde\varphi\big]=0.
\end{aligned}
\end{equation}
Moreover, tangential rotations preserve~$|x-y|$, the halfplanes~$\R^n_\pm$ and~$u_s$ by~\eqref{REARR:in}. Therefore~\eqref{decomp:reflection-invariance} makes~\eqref{decomp:odd} valid in every rotated tangential coordinate system.

Now, using~\eqref{decomp:gradient-orthogonality},
\begin{equation}\label{decomp:directional-norm}
\begin{split}
&\int_{\R^n_+}\left(\sum_{i=1}^{n-1}\theta_i\partial_{x_i}u_s\right)^2\,dx=\sum_{i,j=1}^{n-1}\theta_i\theta_j
\int_{\R^n_+}\partial_{x_i}u_s\,\partial_{x_j}u_s\,dx\\
&\quad=\sum_{i=1}^{n-1}\theta_i^2
\int_{\R^n_+}(\partial_{x_1}u_s)^2\,dx
=\int_{\R^n_+}(\partial_{x_1}u_s)^2\,dx>0.
\end{split}
\end{equation}

Also, thanks to~\eqref{REARR:in} we know that~$u_s(x)=u_s(x-2(x'\cdot\theta)(\theta,0))$.
Thus, for~$i=1,\ldots,n-1$ and~$y_n\ne0$,
\begin{equation}\label{bchkwqgirt34tiybkhaw45678}
\begin{split}
\partial_{x_i}u_s(y)
&=\partial_{x_i}\Big(u_s(y-2(y'\cdot\theta)(\theta,0))\Big)\\
&=\partial_{x_i}u_s(y-2(y'\cdot\theta)(\theta,0))-2\theta_i\sum_{j=1}^{n-1}\theta_j
\partial_{x_j}u_s(y-2(y'\cdot\theta)(\theta,0)).
\end{split}
\end{equation}
We now would like to use this 
with~$y:=x-2(x'\cdot\theta)(\theta,0)$. For this, we observe that
\begin{eqnarray*}
x-2(x'\cdot\theta)(\theta,0) - 2(x'-2(x'\cdot\theta)\theta)\cdot\theta (\theta,0)=x-2(x'\cdot\theta)(\theta,0)
+2(x'\cdot\theta)(\theta,0)=x.
\end{eqnarray*}
As a result, we obtain from~\eqref{bchkwqgirt34tiybkhaw45678} that
\begin{equation}\label{decomp:gradient-reflection}
\partial_{x_i}u_s(x)
-\partial_{x_i}u_s(x-2(x'\cdot\theta)(\theta,0))
=2\theta_i\sum_{j=1}^{n-1}\theta_j\partial_{x_j}u_s(x).
\end{equation}

Multiplying~\eqref{decomp:gradient-reflection} by~$\theta_i$ and summing over~$i$ gives
\begin{equation*}
\begin{split}
\sum_{i=1}^{n-1}\theta_i
\partial_{x_i}u_s(x-2(x'\cdot\theta)(\theta,0))&=\sum_{i=1}^{n-1}\theta_i\partial_{x_i}u_s(x)
-2\left(\sum_{i=1}^{n-1}\theta_i^2\right)
\sum_{j=1}^{n-1}\theta_j\partial_{x_j}u_s(x)\\
&=-\sum_{i=1}^{n-1}\theta_i\partial_{x_i}u_s(x).
\end{split}
\end{equation*}
Therefore, recalling also~\eqref{decomp:gradient-orthogonality} and~\eqref{decomp:coefficients},
\begin{equation}\label{decomp:directional-pairing}
\begin{aligned}
&\frac12\int_{\R^n_+}
\Big(w(x)-w(x-2(x'\cdot\theta)(\theta,0))\Big)
\sum_{i=1}^{n-1}\theta_i\partial_{x_i}u_s(x)\,dx\\
&\quad=\frac12\sum_{i=1}^{n-1}\theta_i
\int_{\R^n_+}w(x)\Big[
\partial_{x_i}u_s(x)
-(\partial_{x_i}u_s)(x-2(x'\cdot\theta)(\theta,0))\Big]\,dx\\
&\quad=\int_{\R^n_+}w\sum_{i=1}^{n-1}\theta_i\partial_{x_i}u_s\,dx
=\sum_{i=1}^{n-1}\theta_i\int_{\R^n_+}w\,\partial_{x_i}u_s\,dx\\
&\quad=\sum_{i=1}^{n-1}c_i\theta_i
\int_{\R^n_+}(\partial_{x_i}u_s)^2\,dx
=\left(\sum_{i=1}^{n-1}c_i\theta_i\right)
\int_{\R^n_+}(\partial_{x_1}u_s)^2\,dx.
\end{aligned}
\end{equation}

Now, by~\eqref{decomp:reflection-invariance}, the function
\begin{equation}\label{cbjwqkefigu234567uiSDFGH}
w^\star(x):=\frac12\big(w(x)-w(x-2(x'\cdot\theta)(\theta,0))\big)
\end{equation}is
a weak solution of~\eqref{lin_problem} in~$H^s_+$ that is odd under the reflection with normal~$\theta$.

Therefore, by~\eqref{decomp:odd}, applied in tangential coordinates with~$\theta$ as a coordinate direction, we see that
\begin{equation}\label{cbjwqkefigu234567uiSDFGH2}
w^\star=\kappa\,\sum_{i=1}^{n-1}\theta_i\partial_{x_i}u_s.\end{equation}
for some~$\kappa\in\R$.

We now claim that
\begin{equation}\label{cbhewgifwgi345678jhgfew}
\kappa=
\sum_{i=1}^{n-1}c_i\theta_i.
\end{equation} To check this,
we use~\eqref{decomp:directional-norm},
\eqref{decomp:directional-pairing} and~\eqref{cbjwqkefigu234567uiSDFGH} to see that
\begin{eqnarray*}&
&\kappa\int_{\R^n_+}(\partial_{x_1}u_s)^2\,dx
=\kappa\int_{\R^n_+}
\left(\sum_{i=1}^{n-1}\theta_i\partial_{x_i}u_s\right)^2\,dx=\int_{\R^n_+}w^\star
\sum_{i=1}^{n-1}\theta_i\partial_{x_i}u_s\,dx\\
&&\qquad=\frac12\int_{\R^n_+}
\Big(w(x)-w(x-2(x'\cdot\theta)(\theta,0))\Big)
\sum_{i=1}^{n-1}\theta_i\partial_{x_i}u_s(x)\,dx\\
&&\qquad=\left(\sum_{i=1}^{n-1}c_i\theta_i\right)
\int_{\R^n_+}(\partial_{x_1}u_s)^2\,dx.
\end{eqnarray*}
Therefore, the strict positivity in~\eqref{decomp:directional-norm} permits
a cancellation, and
this establishes~\eqref{cbhewgifwgi345678jhgfew}.

As a consequence of~\eqref{cbjwqkefigu234567uiSDFGH},
\eqref{cbjwqkefigu234567uiSDFGH2} and~\eqref{cbhewgifwgi345678jhgfew}, we obtain that
\begin{equation*}
\frac12\big(w(x)-w(x-2(x'\cdot\theta)(\theta,0))\big)
=\left(\sum_{i=1}^{n-1}c_i\theta_i\right)
\sum_{j=1}^{n-1}\theta_j\partial_{x_j}u_s(x).
\end{equation*}
From this and~\eqref{decomp:gradient-reflection}, we infer that
\begin{eqnarray*}
&&w(x)-w(x-2(x'\cdot\theta)(\theta,0))=2
\left(\sum_{i=1}^{n-1}c_i\theta_i\right)
\sum_{j=1}^{n-1}\theta_j\partial_{x_j}u_s(x)\\&&\qquad=
\sum_{i=1}^{n-1}c_i\partial_{x_i}u_s(x)
-\sum_{i=1}^{n-1}c_i\partial_{x_i}u_s(x-2(x'\cdot\theta)(\theta,0)),
\end{eqnarray*} that is,
\begin{equation}\label{decomp:remainder-reflection}
\begin{split}
&w(x)-\sum_{i=1}^{n-1}c_i\partial_{x_i}u_s(x)\\
&\quad=w(x-2(x'\cdot\theta)(\theta,0))-\sum_{i=1}^{n-1}c_i
\partial_{x_i}u_s(x-2(x'\cdot\theta)(\theta,0)).
\end{split}
\end{equation}
Accordingly, by the arbitrarity of~$\theta$ we deduce that~$w-\displaystyle\sum_{i=1}^{n-1}c_i\partial_{x_i}u_s$ is radial in~$x'$.

Moreover, by~\eqref{W:LS2} and~\eqref{UYUYi8}, for almost every~$(r,x_n)\in(0,+\infty)\times\R$,
\begin{equation*}
\begin{split}
&\frac1{|\mathbb S^{n-2}|}\int_{\mathbb S^{n-2}}
\left(w(r\omega,x_n)
-\sum_{i=1}^{n-1}c_i\partial_{x_i}u_s(r\omega,x_n)\right)\,d\omega\\
&\quad=\frac1{|\mathbb S^{n-2}|}\int_{\mathbb S^{n-2}}w(r\omega,x_n)\,d\omega
=f_0(r,x_n).
\end{split}
\end{equation*}
Therefore, the radial remainder identified by~\eqref{decomp:remainder-reflection} equals~$f_0$, which proves~\eqref{W:LS}.
\end{proof}

We will also need the following convergence result:

\begin{lemma}\label{lem:weighted-compactness}
Let~$s_j\nearrow1$
and let~$u_{s_j}\in\mathcal M_{s_j}$ be positive minimizers centered as in~\eqref{REARR:in}. Let~$U$ be the classical minimizer given in Lemma~\ref{nustonu1}. 

Let~$f_j\in H^{s_j}_+$. Suppose that
\begin{equation}\label{cjwry34u345678979605878945cdsjh1}
\sup_{j\in\N}\|f_j\|_{s_j,+}<+\infty\end{equation}
and
\begin{equation}\label{cjwry34u345678979605878945cdsjh2}
f_j\to0\text{ in $L^2(B_R\cap\R^n_+)$, for every $R>0$}.\end{equation}

Then
\begin{equation}\label{nondeg:weighted}
\lim_{j\to+\infty}\int_{\R^n_+}
\big(u_{s_j}^{p-1}+U^{p-1}\big)f_j^2\,dx=0.
\end{equation}
\end{lemma}

\begin{proof} For~$j$ sufficiently large, \eqref{ground:sobolev} applies with any fixed exponent in~$(p+1,2_1^*)$. Then, interpolation with~$L^2(\R^n_+)$ and~\eqref{cjwry34u345678979605878945cdsjh1} yield
\begin{equation}\label{nondeg:sobolev}
\|f_j\|_{L^{p+1}(\R^n_+)}\le C\|f_j\|_{s_j,+}\le C,
\end{equation} up to renaming~$C>0$, which is independent of~$j$. 

Also, Lemma~\ref{3e4rfOSJL:LEM} and~\eqref{MCOSM:AS} give uniform bounds for~$u_{s_j}$ and~$U$ in~$L^\infty(\R^n_+)$. Hence, by H\"older's inequality and~\eqref{nondeg:sobolev}, for every~$R>0$,
\begin{equation*}
\int_{\R^n_+}\big(u_{s_j}^{p-1}+U^{p-1}\big)f_j^2\,dx
\le C\|f_j\|_{L^2(B_R\cap\R^n_+)}^2
+C\left(\|u_{s_j}\|_{L^{p+1}(\R^n_+\setminus B_R)}^{p-1}
+\|U\|_{L^{p+1}(\R^n_+\setminus B_R)}^{p-1}\right).
\end{equation*}
The convergences in~\eqref{MCOSM:AS} and~\eqref{cjwry34u345678979605878945cdsjh2} therefore yield that
\begin{equation*}
\limsup_{j\to+\infty}\int_{\R^n_+}
\big(u_{s_j}^{p-1}+U^{p-1}\big)f_j^2\,dx
\le 2C\|U\|_{L^{p+1}(\R^n_+\setminus B_R)}^{p-1}.
\end{equation*}
Since~$U\in L^{p+1}(\R^n_+)$, sending~$R\to+\infty$ proves~\eqref{nondeg:weighted}, as desired.
\end{proof}

We can now prove the main result of this section:

\begin{proof}[Proof of Theorem~\ref{th:nondegeneracy}]
By Theorem~\ref{ADDITIO:PR}, we can take the ground states positive and centered as in~\eqref{REARR:in}. 

Furthermore, by~\eqref{W:LS} and~\eqref{0qoejld23:TSBqP},
we know that
$$ f_0(|x'|,x_n)=w(x',x_n)-\sum_{i=1}^{n-1}c_i\partial_{x_i}u_s$$
solves~\eqref{lin_problem} weakly, by linearity.
Hence, thanks to Proposition~\ref{decompw}, it suffices to prove that there are no nonzero weak solutions of~\eqref{lin_problem} which are radial in~$x'$.

Suppose now, by contradiction, that nonzero radial weak solutions of~\eqref{lin_problem} exist for some sequence~$s_j\nearrow1$. We can normalize
these functions and thus consider the case~$w_{s_j}\in H^{s_j}_+$, radial in~$x'$, with
\begin{equation}\label{nondeg:normalization}
\|w_{s_j}\|_{s_j,+}=1.
\end{equation}

We use~\eqref{WeakForm_eq} with~$s=s_j$ and~$w=w_{s_j}$.
Therefore, Lemma~\ref{lem:energy-compactness}
(which can be used here thanks to~\eqref{nondeg:normalization}), applied with~$G_j=G=\R^n_+$ and~$f_j=w_{s_j}$, produces a function~$w\in H^1(\R^n_+)$ such that, up to a subsequence,
\begin{equation}\label{eq:convw}
w_{s_j}\to w\quad\text{in }L^2(B_R\cap\R^n_+)
\quad\text{for every }R>0.
\end{equation}
We point out that, for all~$i\in\{1,\dots,n-1\}$,
\begin{equation}\label{eq:radiality}
w(x-2x_i e_i)=w(x)
.\end{equation}

Now, by~\eqref{MCOSM:AS} in Lemma~\ref{nustonu1}
we infer that~$u_{s_j}\to U$ almost everywhere.
Thus, in light of Lemma~\ref{3e4rfOSJL:LEM}, the
dominated convergence theorem and~\eqref{eq:convw}, we obtain that, for each~$\varphi\in C_0^\infty(\R^n)$, 
\begin{equation}\label{nondeg:potential-limit}
\begin{split}
&\lim_{j\to+\infty}\left|\int_{\R^n_+}
\big(u_{s_j}^{p-1}w_{s_j}-U^{p-1}w\big)\varphi\,dx\right|\\
&\quad\le \lim_{j\to+\infty}
\|(u_{s_j}^{p-1}-U^{p-1})\varphi\|_{L^2(\R^n_+)}
\|w_{s_j}\|_{L^2(\R^n_+)}\\&\qquad\qquad+
\|U^{p-1}\varphi\|_{L^2(\R^n_+)}
\|w_{s_j}-w\|_{L^2(\operatorname{supp}\varphi\cap\R^n_+)}=
0.
\end{split}
\end{equation}

We now apply~\eqref{compact:mixed} with~$f_j=w_{s_j}$ and use the convergences~\eqref{nondeg:potential-limit} and~\eqref{nu:convergence} to pass equation~\eqref{WeakForm_eq} to the limit as~$j\to+\infty$. In this way, for every~$\varphi\in C_0^\infty(\R^n)$, 
\begin{equation}\label{eq:convw1}
\int_{\R^n_+}(\nabla w\cdot\nabla\varphi+w\varphi)\,dx
=p\nu_1\int_{\R^n_+}U^{p-1}w\varphi\,dx,
\end{equation}
which is the weak Neumann linearized equation at~$U$.

Moreover, the even reflection~$x\mapsto w(x',|x_n|)$ belongs to~$H^1(\R^n)$. Hence, testing~\eqref{eq:convw1} with~$x\mapsto\varphi(x',x_n)+\varphi(x',-x_n)$, for~$\varphi\in C_0^\infty(\R^n)$, and changing variables~$x_n\mapsto-x_n$ in the integrals containing~$\varphi(x',-x_n)$, we see that~$ w(x',|x_n|)$ lies in the kernel of
\begin{equation*}
-\Delta+1-p\nu_1U(x',|x_n|)^{p-1}.
\end{equation*}
Therefore, by~\cite[Theorem~1.2]{FKNOTES}, applied with~$q=p+1$ to the even reflection of~$\nu_1^{1/(p-1)}U$ with the normalization in Remark~\ref{relationnu}, the even reflection of~$w$ is a linear combination of all~$n$ partial derivatives of the even reflection of~$U$. Also, the even reflection of~$w$ and the tangential derivatives of the even reflection of~$U$ are even in~$x_n$, whereas the normal derivative of the even reflection of~$U$ is odd in~$x_n$. Consequently, the normal derivative term vanishes and, as a result,
\begin{equation}\label{nondeg:local-kernel}
w=\sum_{i=1}^{n-1}c_i\partial_{x_i}U.
\end{equation}

By Lemma~\ref{nustonu1}, we know
that~$\partial_{x_i}U$ is odd in~$x_i$ and even in the other tangential coordinates. Hence, subtracting the reflection of~\eqref{nondeg:local-kernel} across~$\{x_i=0\}$ and using~\eqref{eq:radiality}, we conclude that
\begin{equation}\label{nondeg:parity}
0=w(x)-w(x-2x_i e_i)=2c_i\partial_{x_i}U(x).
\end{equation}
Since~$U>0$ and~$U\in L^2(\R^n_+)$, this gives that~$c_1=\dots=c_{n-1}=0$, and accordingly
\begin{equation}\label{nondeg:zero}
w\equiv0.
\end{equation}

Thanks to~\eqref{nondeg:normalization},~\eqref{eq:convw} and~\eqref{nondeg:zero}, we can apply
Lemma~\ref{lem:weighted-compactness} with~$f_j=w_{s_j}$. Thus, taking~$\varphi=w_{s_j}$ in~\eqref{WeakForm_eq} and using~\eqref{nondeg:weighted} and~\eqref{nu:convergence}, we conclude that
\begin{equation}\label{nondeg:contradiction}
1=\lim_{j\to+\infty}\|w_{s_j}\|_{s_j,+}^2
=\lim_{j\to+\infty}p\nu_{s_j}\int_{\R^n_+}u_{s_j}^{p-1}w_{s_j}^2\,dx
=0,
\end{equation}which provides the desired contradiction.
\end{proof}

\section{Uniqueness and proof of Theorem~\ref{th:uniqueness}}\label{uniqueness:section}
 
We now deal with the uniqueness statement.

\begin{proof}[Proof of Theorem~\ref{th:uniqueness}]
By Theorem~\ref{ADDITIO:PR} and the normalization~\eqref{problem_nu}, it suffices to prove uniqueness for positive normalized minimizers centered as in~\eqref{REARR:in}. Suppose, by contradiction, that there are
a sequence~$s_j\nearrow1$ and two distinct minimizers~$u_{s_j}$ and~$v_{s_j}$ for every~$j$. By Lemma~\ref{nustonu1}, both~$u_{s_j}$ and~$v_{s_j}$ converge to~$U$ in~$L^{p+1}(\R^n_+)$ and~$\nu_{s_j}\to\nu_1$.

Let also
\begin{equation*}
d_s:=u_s-v_s,\qquad w_j:=\frac{d_{s_j}}{\|d_{s_j}\|_{s_j,+}}\qquad{\mbox{and}}\qquad
a_s:=p\int_0^1\big((1-\tau)v_s+\tau u_s\big)^{p-1}\,d\tau.
\end{equation*}
Since~$\|w_j\|_{s_j,+}=1$, we deduce from~\eqref{nondeg:sobolev} that$$\sup_{j\in\N}\|w_j\|_{L^{p+1}(\R^n_+)}<+\infty.$$
Also, we subtract the two equations~\eqref{problem_nu} for~$u_{s_j}$ and~$v_{s_j}$, finding that, for every~$\phi\in H^{s_j}_+$,
\begin{equation}\label{OA:diuaL}
(w_j,\phi)_{s_j,+}=\nu_{s_j}\int_{\R^n_+}a_{s_j}w_j\phi\,dx.
\end{equation}

Also, for all~$a,b\ge0$,
$$|a^{p-1}-b^{p-1}|\le\begin{cases}|a-b|^{p-1}&{\mbox{if~$1<p\le2$,}}\\
C(a^{p-2}+b^{p-2})|a-b|
&{\mbox{if~$p>2$,}}\end{cases}$$
and therefore, by H\"older's inequality with~$m:=\frac{p+1}{p-1}$,
\begin{equation}\label{EWQ:0213er}\begin{split}&
\lim_{j\to+\infty}\|a_{s_j}-pU^{p-1}\|_{L^m(\R^n_+)}
\\&\quad\le C\lim_{j\to+\infty}\big(\|u_{s_j}-U\|_{L^{p+1}(\R^n_+)}
+\|v_{s_j}-U\|_{L^{p+1}(\R^n_+)}\big)^{\min\{1,p-1\}}
=0.\end{split}\end{equation}

Also, for~$G_j=G=\R^n_+$, the energy~$E_j$ in~\eqref{compact:forms} satisfies
\begin{equation*}
E_j(w_j)=\|w_j\|_{s_j,+}^2=1.
\end{equation*}
Therefore, Lemma~\ref{lem:energy-compactness} applies with~$f_j=w_j$.
Consequently, up to a subsequence, there exists~$w\in H^1(\R^n_+)$ such that~$w_j\to w$ in~$L^2(B_R\cap\R^n_+)$ for every~$R>0$.
Also, the limit~$w$ retains the tangential symmetry in~\eqref{REARR:in}.

Thus, for all~$\Phi\in C^\infty_0(\R^n)$, we deduce from~\eqref{EWQ:0213er} that
\begin{equation*}
\lim_{j\to+\infty}
\left|\int_{\R^n_+}(a_{s_j}-pU^{p-1})w_j\Phi\,dx\right|
\le\lim_{j\to+\infty}\|a_{s_j}-pU^{p-1}\|_{L^m(\R^n_+)}
\|w_j\|_{L^{p+1}(\R^n_+)}\|\Phi\|_{L^{p+1}(\R^n_+)}
=0.\end{equation*}
This observation allows us to replace~$a_{s_j}$ by~$pU^{p-1}$ when passing to the limit in~\eqref{OA:diuaL}. 

Therefore, using the local~$L^2$ convergence of~$w_j$,~\eqref{compact:mixed} and~\eqref{nu:convergence}, we obtain the classical weak equation~\eqref{eq:convw1} for~$w$. Also, the tangential symmetry of~$w$ gives~\eqref{eq:radiality}. Hence the kernel characterization~\eqref{nondeg:local-kernel} and the reflection identity~\eqref{nondeg:parity}  apply and yield that
\begin{equation}\label{WvanID}
w\equiv0.\end{equation}
By~\eqref{WvanID}, we know that~$w_j\to0$ in~$L^2(B_R\cap\R^n_+)$ for every~$R>0$. Thus, since~$\|w_j\|_{s_j,+}=1$, Lemma~\ref{lem:weighted-compactness} with~$f_j=w_j$ gives
\begin{equation}\label{IL:lim}\lim_{j\to+\infty}\int_{\R^n_+}U^{p-1}w_j^2\,dx=0.\end{equation}

Additionally, testing~\eqref{OA:diuaL} with~$\phi=w_j$, 
recalling~\eqref{EWQ:0213er} and~\eqref{IL:lim},
and passing to the limit as~$j\to+\infty$, we get that
\begin{equation*}1= \lim_{j\to+\infty}\nu_{s_j}\int_{\R^n_+}a_{s_j}w_j^2\,dx
=p\lim_{j\to+\infty}\nu_{s_j}\int_{\R^n_+}U^{p-1}w_j^2\,dx=0.\end{equation*}
This contradiction proves that~$u_{s_j}=v_{s_j}$ for large~$j$,
which gives the desired uniqueness result.
\end{proof}

\section*{Acknowledgements}
SD has been supported by the Australian Future Fellowship~FT230100333,
MM has been supported by the grant PID2023-149451NA-I00 of MCIN/AEI/10.13039/ 501100011033/FEDER, UE,
and EV has been supported by the Australian Laureate Fellowship~FL190100081.

\vfill

\end{document}